\documentclass[journal=jctcce,manuscript=article,layout=tradition]{achemso}

\usepackage{chemformula} 
\usepackage[T1]{fontenc} 
\usepackage{amsmath,amsthm,amssymb}
\usepackage{dsfont}
\usepackage{graphicx}
\usepackage{dcolumn}
\usepackage{bm}
\usepackage{color}
\usepackage{epsfig}
\usepackage{multirow}
\usepackage{mathrsfs}
\usepackage{hyperref}
\usepackage{cleveref}
\usepackage{epstopdf}
\usepackage{subfigure}
\usepackage{autobreak}
\usepackage[utf8]{inputenc}
\usepackage[T1]{fontenc}
\usepackage{mathptmx}
\usepackage{etoolbox}

\newcommand{\V}[1]{\boldsymbol{#1}} 
\newcommand{\T}[1]{\text{#1}} 
\newcommand{\M}[1]{\boldsymbol{#1}} 
\renewcommand{\d}[1]{\delta#1} 
\newcommand{\abs}[1]{\left|#1\right|} 

\newcommand{\grad}{\M{\nabla}} 

\renewcommand{\i}{\mathrm i} 
\newcommand{\eps}{\epsilon}
\newcommand{\mrm}[1]{\mathrm{#1}} 
\newcommand{\te}[1]{\text{#1}}

\author{Xuanzhao Gao}
\altaffiliation{These authors contributed equally}
\affiliation{Thrust of Advanced Materials, and Guangzhou Municipal Key Laboratory of Materials Informatics, The Hong Kong University of Science and Technology (Guangzhou), Guangzhou, China}
\alsoaffiliation[3]{Department of Physics, The Hong Kong University of Science and Technology, Hong Kong SAR, China}
\author{Qi Zhou}
\altaffiliation{These authors contributed equally}
\affiliation{School of Mathematical Sciences, MOE-LSC and CMA-Shanghai, Shanghai Jiao Tong University, Shanghai, China}
\author{Zecheng Gan}
\affiliation{Thrust of Advanced Materials, and Guangzhou Municipal Key Laboratory of Materials Informatics, The Hong Kong University of Science and Technology (Guangzhou), Guangzhou, China}
\alsoaffiliation{Department of Mathematics, The Hong Kong University of Science and Technology, Hong Kong SAR, China}
\email{zechenggan@ust.hk}
\author{Jiuyang Liang}
\affiliation{School of Mathematical Sciences, MOE-LSC and CMA-Shanghai, Shanghai Jiao Tong University, Shanghai, China}
\alsoaffiliation{Center for Computational Mathematics, Flatiron Institute, Simons Foundation, New York, USA}
\email{liangjiuyang@sjtu.edu.cn}

\title[]{Accurate Error Estimates and Optimal Parameter Selection in Ewald Summation for Dielectrically Confined Coulomb Systems}

\abbreviations{IR,NMR,UV}
\keywords{American Chemical Society, \LaTeX}

\begin{document}

\begin{abstract}
Dielectrically confined Coulomb systems are widely employed in molecular dynamics (MD) simulations. Despite extensive efforts in developing efficient and accurate algorithms for these systems, rigorous and accurate error estimates, which are crucial for optimal parameter selection for simulations, is still lacking. In this work, we present a rigorous error analysis in Ewald summation for electrostatic interactions in systems with two dielectric planar interfaces, where the polarization contribution is modeled by an infinitely reflected image charge series. Accurate error estimate is provided for the truncation error of image charge series, as well as decay rates of energy and force correction terms, as functions of system parameters such as vacuum layer thickness, dielectric contrasts, and image truncation levels. Extensive numerical tests conducted across several prototypical parameter settings validate our theoretical predictions. Additionally, our analysis elucidates the non-monotonic error convergence behavior observed in previous numerical studies. Finally, we provide an optimal parameter selection strategy derived from our theoretical insights, offering practical guidance for efficient and accurate MD simulations of dielectric-confined systems.


\end{abstract}


\section{Introduction}

Quasi-two-dimensional (quasi-2D) systems~\cite{mazars2011long}, which are macroscopic in $xy$ dimensions but confined in $z$, have attracted significant attention across various scientific and engineering fields.
Due to confinement effects, the collective behavior of atoms and molecules can differ substantially from that in bulk materials, resulting in intriguing phenomena such as the formation of graphene~\cite{novoselov2004electric}, metal dichalcogenide monolayers~\cite{kumar2012tunable}, and colloidal monolayers~\cite{mangold2003phase}. 
Among quasi-2D systems, charged particles confined by materials with different dielectric constants are of particular interest. 
The so-called \emph{dielectric confinement effect} can further lead to inhomogeneous screening, as compared to bulk systems, or even broken symmetries near interfaces~\cite{C0NR00698J,wang2016inhomogeneous,gao2024broken}.
These effects are crucial in a wide range of applications, including the behavior of water, electrolytes, and ionic liquids confined within thin films~\cite{raviv2001fluidity}, ion transport in nanochannels~\cite{zhu2019ion}, and the self-assembly of colloidal and polymer monolayers~\cite{kim2017imaging}.

Over the past few decades, efficient and accurate simulations of dielectric-confined quasi-2D systems has remained a significant challenge. 
A common approach is by combining the image charge method (ICM)~\cite{jackson1999classical,frenkel2023understanding}, 
where the polarization effects induced by dielectric planar interfaces are modeled through an infinite series of reflected image charges in homogeneous space.
Specifically, this approach truncates the infinitely reflected image charge series along the $z$-direction to a finite number of reflection layers. 
This reduction transforms the original dielectrically inhomogeneous system into a homogeneous one, albeit with an expanded dimension in $z$.
The remaining computational tasks can be effectively handled using the Ewald2D summation method~\cite{parry1975electrostatic,zhonghanhu2014JCTC}. While conceptually straightforward, this ICM-Ewald2D approach is computationally inefficient due to its $O(N^2)$ complexity.
Several alternative approaches have been developed over recent years, including the fast Fourier-Chebyshev spectral method~\cite{maxian2021fast,gao2024fast}, the harmonic surface mapping method~\cite{liang2020harmonic,liang2022hsma}, the random batch Ewald method~\cite{gan2024random}, and boundary element methods~\cite{barros2014efficient,nguyen2019incorporating}.
When combined with the fast multipole method (FMM)~\cite{greengard1987fast,kohnke2020gpu}, the fast Fourier-Chebyshev transform~\cite{trefethen2000spectral}, or random batch importance sampling~\cite{jin2021random}, these methods can achieve computational complexities of $ O(N\log N)$ or even $ O(N)$.

Among existing approaches, a popular technique involves first approximating the exact Ewald2D summation using 3D Ewald summation. 
Then, to eliminate unwanted contributions from periodic replicas along the $z$-axis, correction terms such as the Yeh-Berkowitz (YB) correction~\cite{yeh1999ewald,dos2015electrolytes} and the electric layer correction (ELC)~\cite{arnold2002electrostatics,tyagi2008electrostatic} have been introduced.
Its advantages include: (a) the use of 3D fast Fourier transform (FFT) to achieve $ O(N\log N)$ computational complexity~\cite{yuan2021particle,huang2024pmc}; 
(b) the straightforward incorporation of polarization contributions via image charges, as demonstrated in methods such as ICM-ELC~\cite{tyagi2008electrostatic}, ICM-Ewald3D~\cite{dos2015electrolytes}, ICM-PPPM~\cite{yuan2021particle}, and RBE2D~\cite{gan2024random}; 
and (c) the minimal modifications needed to integrate these quasi-2D electrostatic solvers into mainstream software~\cite{ABRAHAM201519,thompson2021lammps}, where 3D periodic solvers have been extensively optimized.

The aforementioned methods have been widely employed in molecular dynamics (MD) simulations of dielectric-confined systems. It is noteworthy that, for the homogeneous case (i.e., without dielectric mismatch at planar interfaces), rigorous error bounds have been established by Hu and his coworkers~\cite{zhonghanhu2014JCTC,pan2014rigorous}. 
However, a rigorous and accurate error analysis of the correction terms,  particularly concerning their dependence on numerical parameters such as image truncation levels, dielectric contrasts, and vacuum layer thickness, remains elusive. 
Consequently, a systematic framework for optimal parameter selection in MD simulations of dielectric-confined systems is still absent, limiting both accuracy control and practical efficiency.

In fact, the absence of comprehensive error analysis has led to seemingly conflicting conclusions derived solely from numerical tests.
For instance,  dos Santos and Levin~\cite{dos2015electrolytes} suggest that approximately 50 image reflection layers are necessary to achieve satisfactory accuracy, 
whereas Yuan {\it et al.}~\cite{yuan2021particle} argue that only 5 layers suffice. 
Moreover, it has been reported that an excessive number of image reflection layers can lead to counterintuitive losses in accuracy~\cite{yuan2021particle}. This highlights the need for rigorous error analysis to guide optimal parameter selection strategies, thereby addressing the challenges in achieving accurate and efficient MD simulations of dielectric-confined systems.

In this work, we present rigorous and accurate error estimations in Ewald summation for electrostatics of dielectric-confined planar systems. 
First, the truncation error of image charge series, together with the approximation of Ewald2D by 3D Ewald summation formula, is analyzed theoretically and validated by extensive numerical tests. 
Additionally, our theoretical error estimates establish a direct connection to the YB and ELC correction terms, which extends Hu's work~\cite{zhonghanhu2014JCTC,pan2014rigorous} to inhomogeneous cases. 
This development further leads to a systematic parameter optimization scheme for existing ICM-based approaches, such as the ICM-Ewald3D and ICM-PPPM methods. This is particularly advantageous for MD simulations of \emph{strongly-confined} systems.




The subsequent sections of this paper are organized as follows. 
\Cref{sec:background} presents the electrostatic model for dielectric-confined quasi-2D systems, revisiting the image charge and Ewald summation methods.
\Cref{sec:error} analyzes the truncation error of the image charge series and the error associated with reformulating the Ewald2D summation into a 3D Ewald sum, both of which are carefully validated through numerical tests.
\Cref{sec:parameter} provides a parameter selection strategy based on the derived error estimations.
Finally, \Cref{sec:conclusion} summarizes the findings of this work.

\section{Model and methodologies}\label{sec:background}


\subsection{Dielectric-confined quasi-2D systems}\label{sec:system}
Dielectric-confined quasi-2D systems are three-dimensional models characterized by double periodicity in the $xy$ dimensions and confinement by planar dielectric interfaces along the $z$-axis.
The characteristic length scale in the non-periodic direction is often significantly smaller than that in the periodic directions. 
In what follows, we present the mathematical formulation of the model, along with its analytical solution based on image charges.

Consider two parallel dielectric interfaces positioned at $z=0$ and $z=H$, dividing the entire 3D space $\mathbb{R}^3$ into three distinct layers. 
From top to bottom, these layers are denoted as $\Omega_{\text{u}}$, $\Omega_{\text{c}}$ and $\Omega_{\text{d}}$, respectively. 
The central simulation cell $\Omega\in\Omega_{\text{c}}$ has side lengths~$\bm{L}=(L_x, L_y, H)$ and contains $N$ particles located at $\{\V{r}_i\}_{i=1}^{N}$ with charges $\{q_i\}_{i=1}^{N}$. 
The total charge neutrality condition is assumed, i.e., $\sum_{i=1}^{N}q_i=0$. The dielectric permittivity is a piecewise constant function, defined as:
\begin{equation}\label{eq:sanwich}
    \eps(\V{r})=\left\{
        \begin{array}{ll}
        \epsilon_{\T u}, &\, {\V{r} \in \Omega_{\T u}},\\
        \epsilon_{\mrm{c}}, &\, {\V{r} \in \Omega_{\mrm{c}}},\\
        \epsilon_{\mrm{d}}, & \,{\V{r} \in \Omega_{\mrm{d}}},
    \end{array} \right.
\end{equation}
where~$\eps_{\T u}$,~$\eps_{\T c}$ and~$\eps_{\T d}$ are all positive constants, whose specific values depend on materials properties inside each layer. 

The Green's function $G(\V{r},~\V{r}^\prime)$, which describes the electrostatic response at any target location $\V r\in\mathbb{R}^3$ due to a point source charge located at $\V{r}^\prime\in\Omega$, satisfies Poisson's equation with dielectric interface conditions:
\begin{equation}
    \left\{
    \begin{array}{ll}
        - \grad_{\V{r}} \left[ \eps(\V{r}) \grad_{\V{r}} G(\V{r},~\V{r}^\prime) \right] = \d (\V{r} - \V{r}^\prime) & r \in \mathbb{R}^3 \;, \\
        G(\V{r},~\V{r}^\prime) |_{-} = G(\V{r},~\V{r}^\prime) |_{+} & \text{on}~\partial\Omega_{\T c}\;, \\
        \eps_{\T c} G_{\bm{n}}(\V{r},~\V{r}^\prime) |_{-} = \eps_{\T u}  G_{\bm{n}}(\V{r},~\V{r}^\prime) |_{+} & \text{on}~\partial \Omega_{\T c} \cap \partial \Omega_{\T u}\;, \\
        \eps_{\T c} G_{\bm{n}}(\V{r},~\V{r}^\prime) |_{+} = \eps_{\T {d}} G_{\bm{n}}(\V{r},~\V{r}^\prime) |_{-} & \text{on}~\partial \Omega_{\T c} \cap \partial \Omega_{\T d}\;, \\
        G(\V{r},~\V{r}^\prime) \to 0 & \text{as}~r \to \infty\;,
    \end{array}
    \right.
    \label{eq:Poisson_G}
\end{equation}
where $G_{\bm{n}}$ represents the normal derivative of $G$ at planar interfaces, and the subscripts ``$+$/$-$'' denote exterior and interior limits, respectively. 
Due to the double periodicity in $xy$, the total electrostatic energy $U$ can be expressed as
\begin{equation}\label{eq:U_direct}
    U = \frac{1}{2} \sum_{\V{m}} {}^\prime \sum_{i=1}^{N}\sum_{j=1}^{N} q_i q_j G(\V{r}_i,~\V{r}_j + \V{L}_{\V{m}})\;,
\end{equation}
where~$\V{L}_{\V{m}} := (m_x L_x, m_y L_y,0)$ defines the quasi-2D periodic lattice with $\V{m}=(m_x,m_y) \in \mathbb{Z}^2$, and the prime notation $\sum{}^\prime$ indicates that when $i=j$ and $m_x = m_y = 0$, the Green's function should be modified by
\begin{equation}\label{eq::Green}
    G(\V{r}, \V{r}^\prime) \rightarrow \lim_{\V{r} \to \V{r}^\prime} \left(G(\V{r}, \V{r}^\prime) - \frac{1}{4\pi \eps_{\T c}\abs{\V{r} - \V{r}^\prime}} \right)\;,
\end{equation}
 so that the unwanted singular self-interaction is excluded in the summation.

The computation of the Green's function $G$, 
as defined by Eqs.~\eqref{eq:Poisson_G} and~\eqref{eq::Green} is challenging. 
Among existing methodologies, the ICM is frequently used to represent $G$. Specifically, due to the presence of two dielectric planar interfaces, this approach expresses $G$ as an infinite series of reflected images:
\begin{equation}\label{eq:ICM}
    G(\V r,\V r')=\frac{1}{4 \pi \eps_{\mrm{c}}} \left[ \frac{1}{\abs{\V{r} - \V{r}^\prime}} + \sum_{l = 1}^\infty \left( \frac{\gamma_{+}^{(l)}}{\abs{\V{r} - \V{r}_{+}^{\prime(l)}}} + \frac{\gamma_{-}^{(l)}}{\abs{\V{r} - \V{r}_{-}^{\prime (l)}}} \right) \right]\;,
\end{equation}
where~
$\gamma_{+}^{(l)} = \gamma_{\T d}^{\lceil l/2 \rceil} \gamma_{\T u}^{\lfloor l/2 \rfloor}$, $\gamma_{-}^{(l)} = \gamma_{\T d}^{\lfloor l/2 \rfloor} \gamma_{\T u}^{\lceil l/2 \rceil}$, $\V{r}_{+}^{(l)} = (x, y, z_{+}^{(l)})$, and $\V{r}_{-}^{(l)} = (x, y, z_{-}^{(l)})$ are the scaling factors and positions of the $l$-th level image charges.
Here, the notation $\lceil x\rceil$ $(\lfloor x\rfloor)$ represents the ``ceil'' (``floor'') function. ``$+$'' and ``$-$'' indicate that the images are located in~$\Omega_{\mrm{u}}$ and~$\Omega_{\mrm{d}}$, respectively. $\gamma_{\T u}$ and $\gamma_{\T d}$ are reflection factors for the upper and lower dielectric interfaces, defined as:
\begin{equation}
    \gamma_{\mrm{u}} = \frac{\eps_{\mrm{c}} - \eps_{\mrm{u}}}{\eps_{\mrm{c}} + \eps_{\mrm{u}}}~\quad\text{and}\quad
    \gamma_{\mrm{d}} = \frac{\eps_{\mrm{c}} - \eps_{\mrm{d}}}{\eps_{\mrm{c}} + \eps_{\mrm{d}}}\;.
\end{equation}
Finally, the detailed positions of the $l$-th level images along $z$-axis are given by 
\begin{equation}\label{eq:z_l}
    z_{+}^{(l)} = (-1)^l z + 2\lceil l/2\rceil H \;,\quad z_{-}^{(l)} = (-1)^l z - 2\lfloor l/2\rfloor H \;.
\end{equation}
In practice, since all dielectric constants are positive, it follows that $\abs{\gamma_{\mrm{u}}},\abs{\gamma_{\T d}}\leq 1$. This implies that the infinite image charge series converges and can be truncated at the~$M$-th level. 
Consequently, the original dielectric-confined system can be approximated by a homogeneous system augmented with $M$ additional levels of image charges in $z$, which can be readily solved using standard methods for homogeneous quasi-2D systems.

\subsection{The ICM-Ewald2D summation and its reformulation}\label{sec:ewald}

In this section, we first revisit the widely used ``Ewald splitting'' technique~\cite{ewald1921berechnung}, the Ewald2D summation for homogeneous quasi-2D systems, and provide its extension to dielectric-confined cases, i.e., the evaluation of Eq.~\eqref{eq:U_direct}.
Then we introduce a reformulation of the ICM-Ewald2D summation, enabling its efficient computation accelerated via FFT, thereby reducing the cost to $ O(N \log N)$. 

In the Ewald splitting, the Coulomb kernel is decomposed into a sum of short-range and long-range components:
\begin{equation}\label{eq:ewald_decomposition}
\frac{1}{r}=\frac{\mathrm{erfc}(\alpha r)}{r}+\frac{\mathrm{erf}(\alpha r)}{r} \;,
\end{equation}
where $\alpha>0$ is the splitting factor, the error function $\mathrm{erf}(x)$ is  defined as
\begin{equation}\label{eq:erf}
    \mathrm{erf}(x):=\frac{2}{\sqrt{\pi}}\int_0^x e^{-u^2}du\;
\end{equation}
and the complementary error function $\mathrm{erfc}(x):=1-\mathrm{erf}(x)$. 
The advantage of Ewald splitting is clear:
the short-range component, although singular, decays rapidly and is thus well-suited for real space computation, whereas the long-range component, being smooth, can be efficiently handled in reciprocal space.

By Ewald splitting, it has been shown that the quasi-2D lattice sum for electrostatic energy $U$ (\emph{without} dielectric interfaces) can be decomposed as $ U = U_{\text{real}}+U_{\text{Fourier}}$, where~\cite{parry1975electrostatic,zhonghanhu2014JCTC}
\begin{equation}\label{eq:ewald2d-1}
    U_{\te{real}}=\frac{1}{2} \sum_{i, j=1}^N \sum_{\V m}{}^\prime q_iq_j\frac{\mathrm{erfc}(\alpha\abs{\V r_{ij}+\V{L}_{\bm{m}} })}{\abs{\V r_{ij}+ \V{L}_{\bm{m}}}} \;,
\end{equation}
\begin{equation}\label{eq:ewald2d-2}
\begin{split}
    U_{\te{Fourier}} = & \frac{\pi}{2L_xL_y}\sum_{i, j=1}^N q_iq_j\sum_{\V h\neq \V 0}\frac{e^{\i \V h\cdot \V r_{ij}}}{h}\mathcal{G}_{\alpha}(h,z_{ij}) - \frac{\alpha}{\sqrt{\pi}}\sum_{i=1}^{N}q_i^2+\mathcal J_0\;.
\end{split}
\end{equation}
Here, $ {\V h}=2\pi(n_x/L_x, n_y/L_y, 0)$ denotes the reciprocal lattice vector with $n_x,n_y\in \mathbb{Z}$ and $h=|\V{h}|$. 
The notation $\i = \sqrt{-1}$ represents the imaginary unit, and
the function $\mathcal{G}_{\alpha}(\cdot,\cdot)$ is defined as
\begin{equation}\label{eq:G_hzij}
\mathcal{G}_{\alpha}(h,z):=\left[e^{h z}\mathrm{erfc}\left(\frac{h}{2\alpha}+\alpha z\right)+e^{-h z}\mathrm{erfc}\left(\frac{h}{2\alpha}-\alpha z\right)\right]\;,
\end{equation}
with the $\V 0$-th mode correction 
\begin{equation}\label{eq:ewald2d-j0-homo}
    \mathcal J_0 = \frac{-\pi}{L_xL_y}\sum_{i, j=1}^N q_iq_j\left[z_{ij}\mathrm{erf}\left(\alpha z_{ij}\right)+\frac{1}{\alpha\sqrt{\pi}}e^{-\alpha^2z^2_{ij}}\right]\;.
\end{equation}
Eqs.~\eqref{eq:ewald2d-1}--\eqref{eq:ewald2d-j0-homo} are the well-known exact Ewald2D summation formulas~\cite{parry1975electrostatic,zhonghanhu2014JCTC}.

Combined with the ICM, the Ewald2D summation can be extended to accommodate for dielectric-confined systems.
In the ICM-Ewald2D summation~\cite{gan2024random}, Eqs.~\eqref{eq:ewald2d-1}-\eqref{eq:ewald2d-2} are modified as
\begin{equation}\label{eq::Uc}
    U_{\te{real}}^{\text{c}} =  \frac{1}{2} \sum_{i,j=1}^N \sum_{\bm{m}}{}^\prime \sum_{l=0}^M q_iq_{j \pm}^{(l)} \frac{\te{erfc}\left(\alpha \left|\bm{r}_i - \bm{r}_{j\pm}^{(l)}+\V{L}_{\V{m}}\right|\right)}{\left|\bm{r}_i-\bm{r}_{j\pm}^{(l)}+\V{L}_{\V{m}}\right|}\;,
\end{equation}
\begin{equation}\label{eq::Fourier2D}
    U^{\te{c}}_{\te{Fourier}} =  \frac{\pi}{2L_xL_y}\sum_{i, j=1}^N \sum_{l = 0}^{M} q_iq_{j \pm}^{(l)} \sum_{\V h\neq \V 0}\frac{e^{\i \V h\cdot \V r_{ij}}}{h}\mathcal{G}_{\alpha}(h,z_i - z_{j \pm}^{(l)})  - \frac{\alpha}{\sqrt{\pi}}\sum_{i=1}^{N}q_i^2+\mathcal{J}^{\te{c}}_0\;.
\end{equation}
Here $q_{j\pm}^{(l)}=\gamma_{\pm}^{(l)}q_j$ are the $l$-th layer image charge  strengths (with $l = 0$ terms indicating the original source charges), and the $\V 0$-th mode correction term should be modified accordingly as 
\begin{equation}\label{eq:ewald2d-j0}
\mathcal J_0^c = -\frac{\pi}{L_xL_y}\sum_{i,j=1}^{N}\sum_{l=0}^{M} q_iq_{j\pm}^{(l)}\mathcal{G}_{\alpha}^0(|z_i-z_{j\pm}^{(l)}|),
\end{equation}
with $\mathcal G_{\alpha}^{0}(z) := z\mathrm{erf}\left(\alpha z\right)+(\alpha\sqrt{\pi})^{-1}e^{-\alpha^2z^2}$.
Eqs.~\eqref{eq::Uc}--\eqref{eq:ewald2d-j0} integrate the well-established Ewald2D summation formula~\cite{parry1975electrostatic,zhonghanhu2014JCTC} for homogeneous systems with the ICM representation to account for polarization contributions. 
The resulting ICM-Ewald2D formula effectively performs a quasi-2D lattice summation on a system augmented in $z$ by a factor of $2M+1$. 
By selecting a sufficiently large $M$ and setting the real space and reciprocal space cutoffs as $r_c=s/\alpha$ and $k_c=2s\alpha$, respectively, where $s>0$ is a parameter, 
the error due to cutoffs has been estimated as $\sim O(e^{-s^2}/s^2)$~\cite{gan2024fast}, which decays rapidly as $s$ increases. 
However, the pairwise summation terms (over $i$ and $j$) in Eqs.~\eqref{eq::Fourier2D} and \eqref{eq:ewald2d-j0} still lead to a computational complexity of $O(N^2)$, requiring further acceleration techniques.
A widely used approach for acceleration is to reformulate the Ewald2D summation into a triply-periodic Ewald3D summation.
It is noteworthy that for homogeneous systems, such a reformulation was rigorously established by Pan and Hu~\cite{pan2014rigorous}. 
In what follows, we extend their approach to dielectric-confined systems.

First, one can rewrite the function~$\mathcal{G}_{\alpha}(h,z)$ in Eq.~\eqref{eq:G_hzij} into an integral form:
\begin{equation}\label{eq::Ga-integral}
    \mathcal{G}_{\alpha}(h,z) = \int_{-\infty}^{\infty}\frac{e^{-\frac{h^2}{4\alpha^2}-t^2}}{\frac{h^2}{4\alpha^2}+t^2}e^{2\i \alpha z t} dt\;,
\end{equation}
and analogously,
\begin{equation}
\label{eq::J0-integral}
    \mathcal G_{\alpha}^{0}(z)=-\frac{1}{2\pi\alpha }\int_{-\infty}^{\infty}\frac{e^{-t^2}e^{2\i \alpha zt}-1}{t^2}dt\;.
\end{equation}
Discretizing the integrals in Eqs.~\eqref{eq::Ga-integral} and \eqref{eq::J0-integral} using the trapezoidal rule with mesh size $\pi/(\alpha L_z)$, where $L_z>H$ is a parameter, and substituting the discretized forms into Eq.~\eqref{eq::Fourier2D} gives
\begin{equation}\label{eq::UcFour}
\begin{split}
    U^{\T c}_{\T {Fourier}} = \frac{2\pi}{L_xL_yL_z}\sum_{\bm{k}\neq \bm{0}}\frac{e^{-\frac{k^2}{4\alpha^2}}}{k^2}\rho_{\bm{k}}\bar{\rho}_{\bm{k}}^{M}-\frac{\alpha}{\sqrt{\pi}}\sum_{i=1}^{N}q_i^2+ U_{\T {YB}}^{M}+U_{\T {ELC}}^{M}+U_{\T{Trap}}^{M}\;,
\end{split}
\end{equation}
where $\bm{k}=2\pi(n_x/L_x,n_y/L_y,n_z/L_z)$ denotes the 3D periodic lattice vector, and the structure factors $\rho_{\bm{k}}$ and $\widetilde{\rho}_{\bm{k}}^{M}$ are defined as
\begin{equation}
\begin{split}
    \rho_{\bm{k}} := \sum_{i=1}^{N}q_ie^{\i \bm{k}\cdot \bm{r}_i}\; \quad~\text{and}\quad~
    \widetilde{\rho}_{\bm{k}}^{M} & := \sum_{j=1}^{N}q_j\left[e^{-\i \bm{k}\cdot \bm{r}_i}+\sum_{l=1}^{M}\left(\gamma_{+}^{(l)}e^{-\i \bm{k}\cdot \bm{r}_{j+}^{(l)}}+\gamma_{-}^{(l)}e^{-\i \bm{k}\cdot \bm{r}_{j-}^{(l)}}\right)\right]\;.
\end{split}
\end{equation}
On the RHS of Eq.~\eqref{eq::UcFour}, the first two terms resemble the standard Ewald3D summation (with added \emph{vacuum layer} in $z$), where the second term accounts for the self-energy correction. 
The remaining terms provide the additional components required to correct Ewald3D back to Ewald2D:
\begin{equation}\label{eq:U^M_YB}
U_{\T {YB}}^{M}:=\frac{2\pi}{L_xL_yL_z}\left(\sum_{i=1}^{N}q_{i}z_{i}\right)\sum_{j=1}^{N}q_j\left[z_{j}+\sum_{l=1}^{M}\left(\gamma_{+}^{(l)}z_{j+}^{(l)}+\gamma_{-}^{(l)}z_{j-}^{(l)}\right)\right]\;,
\end{equation}
\begin{equation}\label{eq:U^M_ELC}
U_{\T {ELC}}^{M}:=\frac{2\pi}{L_xL_y}\sum_{i,j=1}^{N}q_iq_j\sum_{\bm{h}\neq \bm{0}} \frac{e^{\i \bm{h}\cdot\bm{r}_{ij}}}{h}\frac{\cosh(hz_{ij})+\mathscr{F}_{\text{ELC}}^M(z_i,z_j)}{1-e^{hL_z}}\;,
\end{equation}
where $\mathscr{F}_{\text{ELC}}^M(z_i,z_j)$ is defined as:
\begin{equation}\label{eq::23}
\mathscr{F}_{\text{ELC}}^M(z_i,z_j):=\sum\limits_{l=1}^{M}\left[\gamma_{+}^{(l)}\cosh(h(z_i-z_{j+}^{(l)}))+\gamma_{-}^{(l)}\cosh(h(z_i-z_{j-}^{(l)}))\right].
\end{equation}
The terms in Eqs.~\eqref{eq:U^M_YB} and~\eqref{eq:U^M_ELC} correspond to the ICM-YB~\cite{yuan2021particle} and ICM-ELC~\cite{tyagi2008electrostatic} corrections, respectively. 
The remainder term, $U^{M}_{\T {trap}}$, emerges from the error introduced by trapezoidal discretization. 
The integrand in Eq.~\eqref{eq::Ga-integral} contains two simple poles at $t=\pm {\T i}h/(2\alpha)$, allowing for the estimation of discretization error using contour integral techniques~\cite{trefethen2014Rev}. 
Additionally, the integrand in Eq.~\eqref{eq::J0-integral} is smooth, 
ensuring spectral convergence of the discretization. 
By applying an analysis analogous to that of Pan and Hu~\cite{pan2014rigorous}, we obtain $|U_{\text{trap}}^{M}|\sim e^{-\alpha^2(L_z-H)^2}$, which becomes negligible for $L_z \gg H$.  

In practical computations, the first term in Eq.~\eqref{eq::UcFour} can be efficiently calculated using fast algorithms such as the FFT~\cite{yuan2021particle}, the periodic FMM~\cite{pei2023fast}, and the random batch importance sampling~\cite{liang2022improved}, achieving computational complexities of $O(N\log N)$ or $O(N)$. 
Note that the ICM-YB correction term $U_{\T {YB}}^{M}$ can be directly computed with a cost of $O(N)$, and the remainder term $U_{\text{trap}}^M$ can be eliminated by appropriately choosing $L_z$.
This parameter selection strategy is employed in the recently proposed ICM-Ewald3D~\cite{dos2015electrolytes} and ICM-PPPM~\cite{yuan2021particle} methods. 
However, it is important to note that this approach overlooks the influence of the ICM-ELC correction term $U_{\T {ELC}}^{M}$, which may introduce significant errors. 
Additionally, a rigorous estimate of the image truncation error is also absent in existing works. 
In this study, we address and unify both sources of error.





\section{Error analysis and numerical validations}\label{sec:error}

In this section, we present a comprehensive error analysis for both energy and force calculations, addressing the truncation error of the image charge series and the errors arising from the reformulation of the ICM-Ewald2D summation. 
We also perform extensive numerical experiments to validate our error analysis.
We will focus on force-related results in the main text, while energy-related findings are provided in Section 1 of the Supporting Information (SI).
It is also important to note that in practical computations, the use of FFT can introduce additional errors due to particle spreading onto the uniform grid and the finite resolution of the grid. 
Since these error sources have been thoroughly analyzed in the literature~\cite{deserno1998mesh,wang2012numerical,liang2023error,wang2016multiple,barnett2019parallel,barnett2021aliasing} and are separable from the error discussed in this work, we refer interested readers to these works for more details.

\subsection{Truncation error of the image charge series}\label{sec:error_image}

In this subsection, we analyze the error introduced by truncating the infinite series of image charges at the $M$-th layer. 
This is achieved by reformulating the summation of the infinite image charge series as a Fourier expansion in the periodic $xy$ dimensions and as a geometric series in $z$. 

First, 
let $f(\bm r)$ be a smooth function for $\bm r\in\mathbb R^3$ which decays at infinity, with its Fourier transform denoted as $\widetilde{f}$, then we define the \emph{doubly-periodization} of $f$ as the following lattice sum:
\begin{equation}
f_{\mathcal M} (\bm r)=\sum_{\boldsymbol{m}\in\mathbb{Z}^2} f(\boldsymbol{r}+\V{L}_{\bm{m}})\;,
\end{equation}
and one has the following Poisson summation formula for $f_{\mathcal M} (\bm r)$:
\begin{equation}\label{eq:poissonsum}
f_{\mathcal M} (\bm r)=\frac{1}{2 \pi L_x L_y} \sum_{\bm{k}_{\bm{\rho}}} \int_{\mathbb{R}} \tilde{f}(\bm{k}_{\bm{\rho}}, k_z) e^{\i \bm{k}_{\bm{\rho}} \cdot \boldsymbol{\rho}} e^{\i k_z z} d k_z,
\end{equation}
where $\bm{r}=(\bm{\rho},z)$ with $\bm{\rho}:=(x,y)$, and $\bm{k}_{\bm{\rho}}=2\pi(n_x/L_x,n_y/L_y)$ with $n_x,n_y\in\mathbb{Z}$. 
Applying the Poisson summation formula Eq.~\eqref{eq:poissonsum} to Eq.~\eqref{eq:U_direct} yields 
\begin{equation}\label{eq::33}
    \begin{split}
        U = & \frac{1}{2}\sum_{i,j=1}^{N}q_iq_j\sum_{\bm{m}}{}^\prime\frac{1}{\left|\bm{r}_{ij}+\V{L}_{\bm{m}}\right|} \\
        & + \frac{\pi}{L_xL_y}\sum_{i,j=1}^{N}q_iq_j\sum_{\bm{h}\neq \bm{0}}\frac{e^{-\i \bm{h}\cdot\bm{r}_{ij}}}{h}  \sum\limits_{l=1}^{\infty}\left[\gamma_+^{(l)}e^{-h\left|z_i-z_{j+}^{(l)}\right|}+\gamma_-^{(l)}e^{-h\left|z_i-z_{j-}^{(l)}\right|}\right] \;,
    \end{split}
\end{equation}
where the first term represents the Coulomb energy of a homogeneous quasi-2D system and the second term accounts for the contribution from the image charges. Note that the contribution of the $\bm{h}=\bm{0}$ mode is contained in the homogeneous term, not relevant in the error analysis here for the image charge series.
Now if the image charge series is truncated at the $M$-th layer, then the discarded truncation error term is given by
\begin{equation}\label{eq:trucation_error}
U_{\text{err}}=\frac{\pi}{L_xL_y}\sum_{i,j=1}^{N}q_iq_j\sum_{\bm{h}\neq \bm{0}}\frac{e^{-\i \bm{h}\cdot\bm{r}_{ij}}}{h}\mathcal{G}_{M}(z_i,z_j)\;,  %
\end{equation}
where 
\begin{equation}\label{eq::rGrrij}
 \begin{split}
        \mathcal{G}_{M}(z_i,z_j)&:=\sum\limits_{l=M+1}^{\infty}\left[\gamma_+^{(l)}e^{-h\left|z_i-z_{j+}^{(l)}\right|}+\gamma_-^{(l)}e^{-h\left|z_i-z_{j-}^{(l)}\right|}\right]\; \\
        &= \sum\limits_{l = M + 1}^{\infty} \left[ \gamma_{\T u}^{\lfloor \frac{l}{2} \rfloor}\gamma_{\T d}^{\lceil \frac{l}{2} \rceil} e^{-h (2 \lceil \frac{l}{2} \rceil H + (-1)^l z_i - z_j)} + \gamma_{\T u}^{\lceil \frac{l}{2} \rceil}\gamma_{\T d}^{\lfloor \frac{l}{2} \rfloor} e^{-h (2 \lfloor \frac{l}{2} \rfloor H - (-1)^l z_i + z_j)} \right]\;.
    \end{split}
\end{equation}

To estimate $U_{\text{err}}$, we use the fact that $\mathcal{G}_{M}(z_i,z_j)$ can be reformulated as a geometric series. By inserting the definitions of $\gamma_{\pm}^{(l)}$ and $z_{j \pm}^{(l)}$ into Eq.~\eqref{eq::rGrrij}, we obtain 

\begin{equation}
    \begin{split}
        \mathcal{G}_{M}(z_i,z_j) = & \sum_{n = \frac{M + 1}{2}}^{\infty} \gamma_{\T u}^n \gamma_{\T d}^n e^{-2nhH} \left[ e^{- h z_{ij}} + e^{h z_{ij}} + \gamma_{\T u} e^{- h (z_i + z_j)} + \gamma_{\T d} e^{- h (2 H  - z_i - z_j)}\right] \\
        = & (\gamma_{\T u} \gamma_{\T d})^{\frac{M + 1}{2}} e^{- (M +1) H h} \frac{\left[e^{-h z_{ij}} + e^{h z_{ij}} + \gamma_{\T u} e^{- h (z_i + z_j)} + \gamma_{\T d} e^{- h (2 H - z_i - z_j)} \right] }{1 - \gamma_{\T u} \gamma_{\T d} e^{-2hH}} \;,
    \end{split}
\end{equation}
for odd $M$, and 
\begin{equation}
    \begin{split}
        \mathcal{G}_M(z_i,z_j) = & \sum_{n = \frac{M}{2}}^{\infty}  \gamma_{\T u}^n \gamma_{\T d}^n e^{- 2nhH} \left[  \gamma_{\T u} \gamma_{\T d} e^{- h (2 H  + z_{ij})} + \gamma_{\T u} \gamma_{\T d}  e^{- h (2 H  - z_{ij})} + \gamma_{\T u} e^{- h (z_i + z_j)} + \gamma_{\T d} e^{- h (2H  - z_i - z_j)} \right] \\
        = & (\gamma_{\T u} \gamma_{\T d})^{\frac{M}{2}} e^{- M H h} \frac{\left[ \gamma_{\T u} \gamma_{\T d} e^{- h (2 H + z_{ij})} + \gamma_{\T u} \gamma_{\T d} e^{- h (2H - z_{ij})} + \gamma_{\T u} e^{ - h (z_i + z_j)} + \gamma_{\T d} e^{-h (2 H - z_i - z_j)} \right]}{1 - \gamma_{\T u} \gamma_{\T d} e^{-2hH}}  \;,
    \end{split}
\end{equation}
for even $M$. In both cases, we obtain the following error bound:
\begin{equation}\label{eq:beta_ij}
    \abs{ \mathcal{G}_{M}(z_i,z_j) } \leq \abs{\gamma_{\T u} \gamma_{\T d} e^{- 2 H h}}^{\lfloor{(M+1)/2}\rfloor} \frac{\abs{\gamma_{\T u}} + \abs{\gamma_{\T d}} + 2 \abs{\gamma_{\T u} \gamma_{\T d}}}{1 - \abs{\gamma_{\T u} \gamma_{\T d}} e^{-2hH}} \leq \frac{4 \abs{\gamma_{\T u} \gamma_{\T d} e^{- 2 H h}}^{\lfloor{(M+1)/2}\rfloor}}{1 - \abs{\gamma_{\T u} \gamma_{\T d}} e^{-2hH}} \;,
\end{equation}
where we use the fact that 
\begin{equation}
    \begin{split}
        & ~ \abs{\gamma_{\T u} \gamma_{\T d}} e^{-h(2 H + z_{ij})} + \abs{\gamma_{\T u} \gamma_{\T d}} e^{ - h (2 H - z_{ij})} + \abs{\gamma_{\T d}} e^{- h (2 H - z_i - z_j)} + \abs{\gamma_{\T u}} e^{- h (z_i + z_j)} \\
        \leq & ~ 2 \abs{\gamma_{\T u} \gamma_{\T d}} +\abs{\gamma_{\T d}} + \abs{\gamma_{\T u}} \leq 4 \;.
    \end{split}
\end{equation}
Substituting Eq.~\eqref{eq:beta_ij} into Eq.~\eqref{eq:trucation_error}, we have

\begin{equation}\label{eq:U_gamma_bound}
    \begin{split}
        \left|U_{\text{err}}\right| & \leq \frac{4 \pi}{L_x L_y}  \sum_{\bm{h}\neq \bm{0}}  \frac{\abs{\gamma_{\T u} \gamma_{\T d} e^{- 2 hH}}^{\lfloor \frac{M + 1}{2} \rfloor}}{h(1 - \abs{\gamma_{\T u} \gamma_{\T d}} e^{- 2hH})}\abs{\sum_{i,j=1}^{N}q_i q_j e^{-\i \bm{h}\cdot\bm{r}_{ij}}} \\
        & \leq \frac{4\pi}{L_x L_y} \frac{C_{q} \abs{\gamma_{\T u} \gamma_{\T d}}^{\lfloor{(M+1)/2}\rfloor}}{1 - \abs{\gamma_{\T u} \gamma_{\T d}} e^{- 4\pi H / \max\{L_x,L_y\}}} \sum_{\bm{h}\neq \bm{0}} \frac{e^{- 2 \lfloor{(M+1)/2}\rfloor hH}}{h}\\
        & \leq \frac{4\pi}{L_x L_y} \frac{C_{q} \abs{\gamma_{\T u} \gamma_{\T d}}^{\lfloor{(M+1)/2}\rfloor}}{1 - \abs{\gamma_{\T u} \gamma_{\T d}} e^{- 4\pi H / \max\{L_x,L_y\}}} \int_{\frac{2\pi}{\max\{L_x,L_y\}}}^{\infty} 2\pi h \frac{e^{- 2 \lfloor{(M+1)/2}\rfloor hH}}{h} \, \mathrm{d}h \\
        & = \frac{8 \pi^2}{L_x L_y} \frac{C_{q} \abs{\gamma_{\T u} \gamma_{\T d}}^{\lfloor{(M+1)/2}\rfloor} e^{-\frac{4 \pi H \lfloor{(M+1)/2}\rfloor}{\max\{L_x,L_y\}}}}{1 - \abs{\gamma_{\T u} \gamma_{\T d}} e^{- 4\pi H / \max\{L_x,L_y\}}},
    \end{split}
\end{equation}
where $C_q$ is the bound of $|\sum_{i,j}q_iq_je^{-\i \bm{h}\cdot\bm{r}_{ij}}|$, which depends on the charge distribution of the system. 
A rough bound is $C_q\leq \sum_{i,j}|q_iq_j|\leq N^2q_{\text{max}}^2$, where $q_{\text{max}}=\max_{i}|q_i|$ denotes the maximum strength of a single point charge. 
This estimate can be further refined based on prior knowledge of the charge distribution. 
For example, under the Debye-H\"uckel (DH) approximation~\cite{levin2002electrostatic,gan2024fast}, the bound can be tightened to $C_q\leq CNq_{\max}^2$, where $C$ is a constant independent of $N$. 
Either way, we have the following rate of convergence for $|U_{\text{err}}|$ in terms of the truncation parameter $M$:
\begin{equation}\label{eq:Uerr}
|U_{\text{err}}|\sim O\left(\abs{\gamma_{\T u} \gamma_{\T d}}^{\lfloor\frac{M+1}{2}\rfloor} e^{-\frac{4 \pi H \lfloor (M + 1)/2\rfloor}{\max\{L_x,L_y\}}}\right)\;.
\end{equation}

Next, we analyze the truncation error of the force exerted on the $i$-th particle. Using the definition $\bm{F}^i_{\T{err}}=-\nabla_{\bm{r}_i} U_{\text{err}}$ and taking the derivative in each direction, we obtain:
\begin{equation}\label{eq::Ferr}
    \left|\bm{F}^i_{\T{err}}\right| = \left|\nabla_{\bm{r}_i} U_{\text{err}}\right|\leq  \abs{\frac{2\sqrt{3}\pi}{L_xL_y} q_i \sum_{j=1}^{N} q_j\sum_{\bm{h}\neq \bm{0}} e^{-\i \bm{h}\cdot\bm{r}_{ij}}\sum\limits_{l = M + 1}^{\infty}\mathcal{G}_{M}(z_i,z_j)}, 
\end{equation}
where we use the identity $\partial_{z_i}\sum_{ij}\mathcal{G}_M(z_{i},z_{j})=2h\sum_{ij}\mathcal{G}_{M}(z_i,z_j)$, and the factor $\sqrt{3}$ accounts for three dimensions of the force field $\bm F$.  
Substituting Eq.~\eqref{eq:beta_ij} into Eq.~\eqref{eq::Ferr}, we have 
\begin{equation}
\begin{split}
 |F^{i}_{\text{err}}| 
 & \leq \frac{8\sqrt{3}\pi}{L_xL_y}  \frac{C_{Q} \abs{\gamma_{\T u} \gamma_{\T d}}^{\lfloor{(M+1)/2}\rfloor}}{1 - \abs{\gamma_{\T u} \gamma_{\T d}} e^{- 4\pi H / \max\{L_x,L_y\}}} \sum_{\bm{h}\neq \bm{0}} e^{- 2 \lfloor{(M+1)/2}\rfloor hH}\\
 & \leq \frac{8\sqrt{3}\pi}{L_xL_y} \frac{ C_Q\abs{\gamma_{\T u} \gamma_{\T d}}^{\lfloor{(M+1)/2}\rfloor}}{1 - \abs{\gamma_{\T u} \gamma_{\T d}} e^{- 4\pi H / \max\{L_x,L_y\}}} \int_{\frac{2\pi}{\max\{L_x,L_y\}}}^{\infty} 2\pi h e^{- 2 \lfloor{(M+1)/2}\rfloor hH} \, \mathrm{d}h \\
& = \frac{16\sqrt{3} \pi^2}{L_xL_y} \frac{C_Q \abs{\gamma_{\T u} \gamma_{\T d}}^{\lfloor{(M+1)/2}\rfloor}}{1 - \abs{\gamma_{\T u} \gamma_{\T d}} e^{- 4\pi H / \max\{L_x,L_y\}}}  \frac{\left[1 + \frac{2\pi \lfloor{(M+1)/2}\rfloor H}{\max\{L_x,L_y\}}\right]}{4 \lfloor{(M+1)/2}\rfloor^2 H^2} e^{- \frac{4\pi \lfloor{(M+1)/2}\rfloor H}{\max\{L_x,L_y\}}}\;,
\end{split}
\end{equation}
where the second inequality is derived from the monotonicity of the exponential function, and $C_Q=|q_{i}\sum_{j=1}^{N}q_je^{-\i\bm{h}\cdot \bm{r}_{ij}}|\leq Nq_{\text{max}}^2$. 
If the DH theory is applied, the prefactor $C_Q$ can be further tightened to $C_Q\leq Cq_{\max}^2$, where $C$ is a constant independent of $N$. 
Similarly, considering the rate of convergence in terms of $M$, we have
\begin{equation}\label{eq::Ferr1}
|\bm{F}^i_{\text{err}}| \sim O\left(\lfloor (M+1)/2\rfloor^{-1}\abs{\gamma_{\T u} \gamma_{\T d}}^{\lfloor\frac{M+1}{2}\rfloor} e^{-\frac{4\pi H \lfloor(M+1)/2\rfloor}{\max\{L_x,L_y\}}}\right)\;.
\end{equation}
Comparing Eq.~\eqref{eq::Ferr1} and Eq.~\eqref{eq:Uerr},
we observed that the error in force converge slightly faster than that of the energy by a factor of $\lfloor (M+1)/2\rfloor^{-1}$. 

To validate our theoretical predictions, we employ the ICM-Ewald2D summation (i.e., Eqs.~\eqref{eq::Uc} and~\eqref{eq::Fourier2D}) to calculate the relative errors in force, $\mathcal{E}_r = \max\limits_{i=1,\ldots,N}\frac{\abs{\bm{F}^i_{\text{err}}}}{\abs{\bm{F}^i}}$, as a function of $M$, the number of truncated image layers.
Without loss of generality, we examine charge-asymmetric systems containing $13$ divalent cations and $26$ monovalent anions, randomly distributed within the simulation cell. 
In all calculations, the dimensions of the simulation cell along the periodic directions are fixed as $L_x = L_y = 10$, while the aspect ratio of the system is adjusted by varying the cell height $H$. 
Recall that the reflection factors for the upper and lower dielectric interfaces are denoted as $\gamma_{\T u}$ and $\gamma_{\T d}$, respectively. 
Finally, unless otherwise specified, we always set the Ewald splitting parameter $s=6$ to ensure that the errors due to Ewald decomposition remain negligible.

\begin{figure}[htbp]
    \centering
    \includegraphics[width=0.98\linewidth]{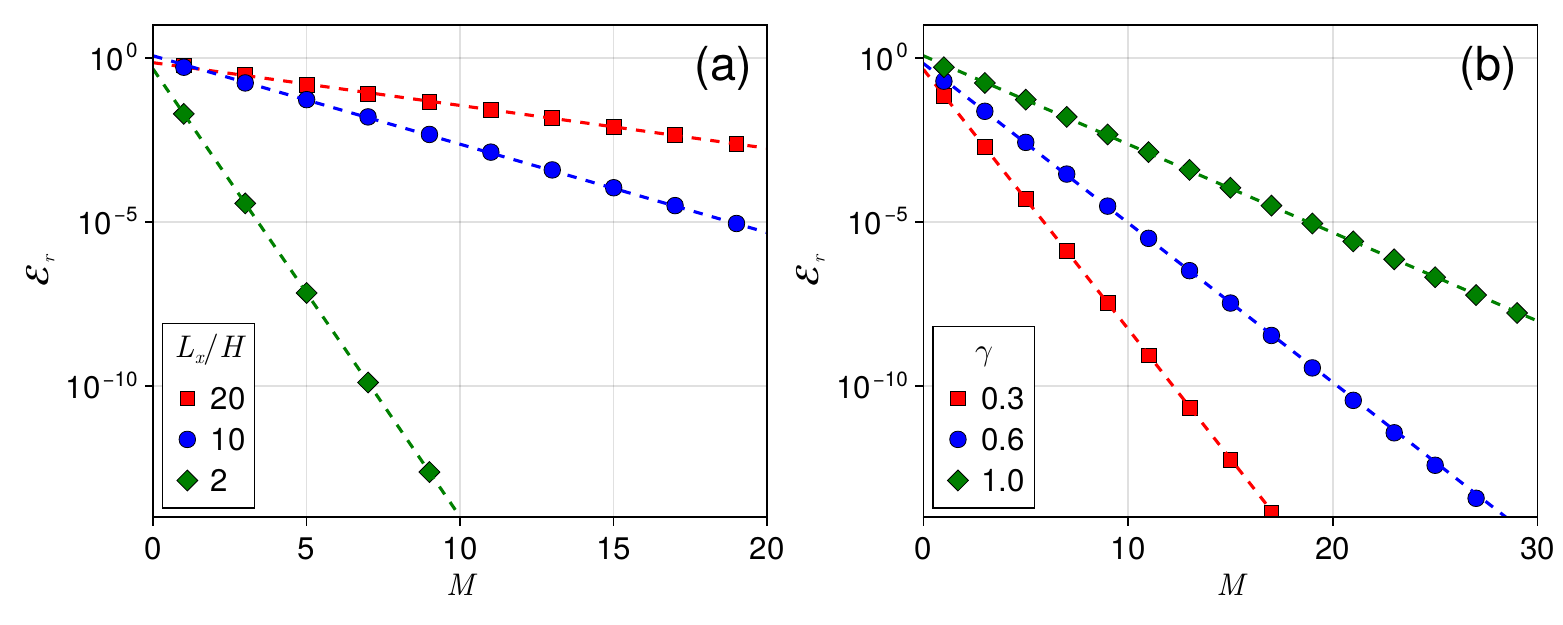}
    \caption{
        Relative errors in force ($\mathcal{E}_r$) as a function of the truncation parameter $M$ for the image charge series. 
        The dashed lines represent the fitted curves with decay rates using our theoretical prediction Eq.~\eqref{eq::Ferr1}. 
        In panel (a), we fix $\gamma_{\T u} = \gamma_{\T d} = 1$ and consider systems with varying heights $H = 0.5$, $1$, and $5$. In panel (b), we fix $H = 1$ while varying $\gamma_{\T u} = \gamma_{\T d} = \gamma$ with values of $0.3, 0.6,$ and $1$. In both panels, we fix $L_x=L_y=10$.
    }
    \label{fig:icm_error_force}
\end{figure}

The numerical results are presented in \Cref{fig:icm_error_force} (a) and (b), where we examine the convergence rate of $\mathcal{E}_r$ under different aspect ratios $L_x / H$ and reflection factors $\gamma=\gamma_{\T u}=\gamma_{\T d}$, respectively.
All results demonstrate that the force errors decay exponentially as $M$ increases. 
Additionally, we observe that the errors decay slower as the aspect ratio $L_x / H$ (\Cref{fig:icm_error_force} (a)) and the reflection factor $\gamma$ (\Cref{fig:icm_error_force} (b)) increases, both are consistent with our theoretical predictions in Eq.~\eqref{eq::Ferr1}. 
Furthermore, to quantitatively validate our theoretical findings, we also plot the theoretical decay rates predicted by Eq.~\eqref{eq::Ferr1} as dashed lines in \Cref{fig:icm_error_force}, showing excellent agreement with our numerical results.
The results for relative errors in energy are documented in Section 1 of the SI, where similar conclusions hold.  


\subsection{Estimations of electrostatic layer correction (ELC) with image charges}\label{sec:error_reform}

Another often-overlooked source of error in existing algorithms stems from the omission of the ELC term in Eq.~\eqref{eq:U^M_ELC} and the discretization error using the trapezoidal rule, which has been discussed in Section~\ref{sec:ewald}. 
In this section, we focus on deriving estimations of the ELC term with image charges.

First, we recall the definition of the ELC term with image charges, denoted as $U_{\text{ELC}}^{M}$, as provided in the energy expressions Eqs.~\eqref{eq:U^M_ELC} and~\eqref{eq::23}:
\begin{equation}\label{eq::UELCm}
U_{\T {ELC}}^{M}:=\frac{2\pi}{L_xL_y}\sum_{i,j=1}^{N}q_iq_j\sum_{\bm{h}\neq \bm{0}} \frac{e^{\i \bm{h}\cdot\bm{r}_{ij}}e^{-hL_z}\left[\cosh(hz_{ij})+\mathscr{F}_{\text{ELC}}^{M}(z_i,z_j)\right]}{h(e^{-hL_z}-1)}\;.
\end{equation}
To derive an estimation for $U_{\T {ELC}}^{M}$, we introduce the following two inequalities. 
1). By basic algebraic manipulations, one obtains:  
\begin{equation}\label{eq::bpud1}
e^{-hL_z}\cosh(hz_{ij})\leq e^{-h(L_z-|z_{ij}|)}\leq e^{-h(L_z-H)}\;,
\end{equation}
and 2).
\begin{equation}\label{eq::bpud2}
\begin{split}
e^{-hL_z}\mathscr{F}_{\text{ELC}}^{M}(z_i,z_j)&=e^{-hL_z}\sum\limits_{l=1}^{M}\left[\gamma_{+}^{(l)}\cosh(h(z_i-z_{j+}^{(l)}))+\gamma_{-}^{(l)}\cosh(h(z_i-z_{j-}^{(l)}))\right]\\
&= \frac{1}{2}\sum\limits_{l=1}^{M}\left[\gamma_{+}^{(l)}e^{-h(L_z-|z_i-z_{j+}^{(l)}|)}+\gamma_{-}^{(l)}e^{-h(L_z-|z_i-z_{j-}^{(l)}|)}\right]+\frac{e^{-hL_z}}{2}(\mathcal{G}_{0}(z_i,z_j)-\mathcal{G}_{M}(z_i,z_j))\\
&\leq \frac{1}{2}\sum_{l=1}^{M}C^{(l)}_{\gamma}e^{-h(L_z-(l+1)H)}+\frac{(|\gamma_{\T u}|+|\gamma_{\T d}| + 2|\gamma_{\T u}\gamma_{\T d}|)e^{-hL_z}}{1-|\gamma_{\T u}\gamma_{\T d}|e^{-2hH}},
\end{split}
\end{equation}
where the factor $C^{(l)}_{\gamma}$ is defined as 
\begin{equation}
C^{(l)}_{\gamma}:=\left|\gamma_{+}^{(l)}\right|+\left|\gamma_{-}^{(l)}\right|=\left| \gamma_{\T d}^{\lceil l/2 \rceil} \gamma_{\T u}^{\lfloor l/2 \rfloor}\right|+\left| \gamma_{\T d}^{\lfloor l/2 \rfloor} \gamma_{\T u}^{\lceil l/2 \rceil}\right|.
\end{equation}
To obtain the above inequality, we use the definitions of $\gamma_{+}^{(l)}$ and $\gamma_{-}^{(l)}$, the bound $\max_{i,j}\{|z_i-z_{j+}^{(l)}|,|z_i-z_{j-}^{(l)}|\}\leq (l+1)H$, and the definition of $\mathcal{G}_{M}(z_i,z_j)$ in Eq.~\eqref{eq::rGrrij} along with its bound given in Eq.~\eqref{eq:beta_ij}. 
Substituting Eqs.~\eqref{eq::bpud1} and \eqref{eq::bpud2} into Eq.~\eqref{eq::UELCm} yields

\begin{equation}
\resizebox{.98\hsize}{!}{$
\begin{split}
\left|U_{\text{ELC}}^{M}\right|&\leq \frac{2\pi}{L_xL_y}\sum_{\bm{h}\neq \bm{0}} \frac{C_q\left[e^{-h(L_z-H)}+\frac{1}{2}\sum\limits_{l=1}^{M}C^{(l)}_{\gamma}e^{-h(L_z-(l+1)H)}+\frac{(|\gamma_{\T u}|+|\gamma_{\T d}| + 2|\gamma_{\T u}\gamma_{\T d}|)e^{-hL_z}}{1-|\gamma_{\T u}\gamma_{\T d}|e^{-\frac{4\pi H}{\max\{L_x,L_y\}}}}\right]}{h(1-e^{-\frac{2\pi L_z}{\max\{L_x,L_y\}}})}\\
&\leq \frac{4\pi^2 C_q}{L_xL_y(1-e^{-\frac{2\pi L_z}{\max\{L_x,L_y\}}})}\int^{\infty}_{\frac{2\pi}{\max\{L_x,L_y\}}} \left[e^{-h(L_z-H)}+\sum\limits_{l=1}^{M}\frac{C^{(l)}_{\gamma}}{2}e^{-h(L_z-(l+1)H)}+\frac{4e^{-hL_z}}{1-e^{-\frac{4\pi H}{\max\{L_x,L_y\}}}}\right]dh\\
&=\frac{4\pi^2 C_q}{L_xL_y(1-e^{-\frac{2\pi L_z}{\max\{L_x,L_y\}}})}\left[\frac{e^{-\frac{2\pi(L_z-H)}{\max\{L_x,L_y\}}}}{L_z-H}+\sum_{l=1}^{M}\frac{C_{\gamma}^{(l)}e^{-\frac{2\pi(L_z-(l+1)H)}{\max\{L_x,L_y\}}}}{2(L_z-(l+1)H)}+\frac{4e^{-\frac{2\pi L_z}{\max\{L_x,L_y\}}}}{L_z(1-e^{-\frac{4\pi H}{\max\{L_x,L_y\}}})} \right],
\end{split}
$}
\end{equation}
where we recall $C_q$ is the bound of $|\sum_{i,j=1}^{N}q_iq_je^{\i \bm{h}\cdot\bm{r}_{ij}}|$, $\left|\gamma_{\T u}\right|\leq 1$, and $\left|\gamma_{\T d}\right|\leq 1$. 
Finally, suppose $L_z>(M+1)H$, we obtain the following estimation for the ELC contribution in energy, $U_{\text{ELC}}^{M}$, as:
\begin{equation}
\label{eq::U_ELC}
\left|U_{\text{ELC}}^{M}\right|\sim O\left(e^{-\frac{2\pi(L_z-H)}{\max\{L_x,L_y\}}}+\sum_{l=1}^{M}C_{\gamma}^{(l)}e^{-\frac{2\pi(L_z-(l+1)H)}{\max\{L_x,L_y\}}}\right).
\end{equation}

The corresponding ELC term in force calculations can be estimated analogously. 
We have
\begin{equation}
\left|\bm{F}_{\text{ELC}}^{M,i}\right|:=|-\nabla_{\bm{r}_i} U_{\text{ELC}}^{M}|\leq \left|\frac{4\sqrt{3}\pi}{L_xL_y}q_i\sum_{j=1}^{N}q_j\sum_{\bm{h}\neq \bm{0}} \frac{e^{\i \bm{h}\cdot\bm{r}_{ij}}e^{-hL_z}\left[\cosh(hz_{ij})+\mathscr{F}_{\text{ELC}}^{M}(z_i,z_j)\right]}{1-e^{-hL_z}}\right|\;.
\end{equation}
Further applying the inequalities Eqs.~\eqref{eq::bpud1} and \eqref{eq::bpud2}, we obtain:

\begin{equation}
\resizebox{0.98\hsize}{!}{$
\begin{split}
\left|\bm{F}_{\text{ELC}}^{M,i}\right|\leq &\frac{4\sqrt{3}\pi}{L_xL_y}\sum_{\bm{h}\neq \bm{0}} \frac{C_Q\left[e^{-h(L_z-H)}+\frac{1}{2}\sum\limits_{l=1}^{M}C^{(l)}_{\gamma}e^{-h(L_z-(l+1)H)}+\frac{(|\gamma_{\T u}|+|\gamma_{\T d}| + 2|\gamma_{\T u}\gamma_{\T d}|)e^{-hL_z}}{1-|\gamma_{\T u}\gamma_{\T d}|e^{-\frac{4\pi H}{\max\{L_x,L_y\}}}}\right]}{1-e^{-\frac{2\pi L_z}{\max\{L_x,L_y\}}}}\\
\leq& \frac{8\sqrt{3}\pi^2 C_Q}{L_xL_y(1-e^{-\frac{2\pi L_z}{\max\{L_x,L_y\}}})}\int_{\frac{2\pi}{\max\{L_x,L_y\}}}^{\infty}h\left[e^{-h(L_z-H)}+\frac{1}{2}\sum\limits_{l=1}^{M}C^{(l)}_{\gamma}e^{-h(L_z-(l+1)H)}+\frac{4e^{-hL_z}}{1-e^{-\frac{4\pi H}{\max\{L_x,L_y\}}}}\right]dh\\
=&\frac{8\sqrt{3}\pi^2 C_Q}{L_xL_y(1-J_3)}\left[\frac{1+J_1}{(L_z-H)^2}e^{-J_1}+\sum_{l=1}^{M}\frac{C_{\gamma}^{(l)}(1+J_2^{(l)})}{2(L_z-(l+1)H)^2}e^{-J_2^{(l)}}+\frac{4(1+J_3)}{L_z^2(1-e^{-\frac{4\pi H}{\max\{L_x,L_y\}}})}e^{-J_3}\right]
\end{split}
$}
\end{equation}
where $C_Q$ is the bound of $|q_i\sum_{j=1}^{N}q_je^{-\i \bm{h}\cdot\bm{r}_{ij}}|$ and the coefficients $J_1$, $J_2^{(l)}$ and $J_3$ are defined via
\begin{equation}
J_1=\frac{2\pi(L_z-H)}{\max\{L_x,L_y\}},\quad\; J_2^{(l)}= \frac{2\pi(L_z-(l+1)H)}{\max\{L_x,L_y\}}\quad \;\text{and} \quad \; J_3=\frac{2\pi L_z}{\max\{L_x,L_y\}}.
\end{equation}
As before, by omitting all the prefactors, we arrive at the following estimation for the ELC contribution in force calculations:
\begin{equation}\label{eq::ForceError}
\left|\bm{F}_{\text{ELC}}^{M,i}\right|\sim O\left(e^{-\frac{2\pi(L_z-H)}{\max\{L_x,L_y\}}}+\sum_{l=1}^{M}C_{\gamma}^{(l)}e^{-\frac{2\pi(L_z-(l+1)H)}{\max\{L_x,L_y\}}}\right).
\end{equation}
Clearly, comparing Eq.~\eqref{eq::ForceError} and Eq.~\eqref{eq::U_ELC}, we find that the ELC contribution behaves asymptotically the same for both energy and force calculations.

\subsection{Leading-order analysis of the ELC term and numerical validations}\label{sec::leadingerr}
The theoretical estimations of the ELC term, as presented in Eqs.~\eqref{eq::ForceError} and~\eqref{eq::U_ELC}, behave differently under different system aspect ratios and reflection factors.
In this section, we conduct a detailed analysis of the leading-order contribution of the ELC term across different system parameter scenarios. 
Our focus will be on the analysis of force calculations, a similar approach can be applied to the energy.

First, by introducing two new dimensionless parameters $g_{\T u}:=\gamma_{\T u}e^{\frac{2\pi H}{\max\{L_x,L_y\}}}$ and $g_{\T d}:=\gamma_{\T d}e^{\frac{2\pi H}{\max\{L_x,L_y\}}}$, Eq.~\eqref{eq::ForceError} can be reformulated as
\begin{equation}
   \label{eq::Error_refo}
   \begin{aligned} \left|\bm{F}_{\text{ELC}}^{M,i}\right|&\sim e^{-\frac{2\pi(L_z-H)}{\max\{L_x,L_y\}}} \left[ 1 + \sum_{l=1}^{M} \left( g_{\T u}^{\lfloor \frac{l}{2} \rfloor} g_{\T d}^{\lceil \frac{l}{2} \rceil} + g_{\T u}^{\lceil \frac{l}{2} \rceil} g_{\T d}^{\lfloor \frac{l}{2} \rfloor} \right) \right]\\
   &=e^{-\frac{2\pi(L_z-H)}{\max\{L_x,L_y\}}}\left[(g_{\T u} + g_{\T d} + 2) \sum_{l = 0}^{\lfloor \frac{M - 1}{2} \rfloor} (g_{\T u}g_{\T d})^{l} + ((-1)^{M}+1) (g_{\T u}g_{\T d})^{\lfloor \frac{M}{2} \rfloor} - 1\right]\,.
   \end{aligned}
\end{equation}
From Eq.~\eqref{eq::Error_refo}, it is clear that the leading-order contribution of the ELC term depends on the magnitude of $|g_{\T u}g_{\T d}|$. 
Specifically, i). if $|g_{\T u}g_{\T d}| > 1$, the ELC leading-order term grows exponentially as $M$ increases:
\begin{equation}\label{eq::lhs_bigger1}
   \abs{\V{F}_{\text{ELC}}^{M,i}} =  \left\{
	\begin{aligned}
		& O \left(2 \abs{g_{\T u} g_{\T d}}^{M/2} e^{-\frac{2\pi (L_z - H)}{\max\{L_x,L_y\}}} \right) , & \text{if $M$ is even,}\\
		& O \left( \left(g_{\T u} + g_{\T d}+2\right)\abs{g_{\T u} g_{\T d}}^{\frac{M-1}{2}} e^{-\frac{2\pi (L_z -H)}{\max\{L_x,L_y\}}} \right) , & \text{if $M$ is odd.}
	\end{aligned}
	\right.
\end{equation}
ii). If $|g_{\T u}g_{\T d}|=1$, the ELC leading-order term grows linearly with $M$:
\begin{equation}
\label{eq::lhs_1}
\abs{\V{F}_{\text{ELC}}^{M,i}} = O \left( \frac{M}{2}\left(g_{\T u} + g_{\T d}+2\right)e^{-\frac{2\pi L_z}{\max\{L_x,L_y\}}} \right)\;.
\end{equation}
iii). If $|g_{\T u}g_{\T d}|<1$, the summation in Eq.~\eqref{eq::Error_refo} converges as $M\rightarrow +\infty$, yielding a uniform estimation independent with $M$:
\begin{equation}
\label{eq::lhs_less1}
\abs{\V{F}_{\text{ELC}}^{M,i}} \sim O \left( \left(g_{\T u}+g_{\T d}+2\right)e^{-\frac{2\pi L_z}{\max\{L_x,L_y\}}} \right).
\end{equation}
The preceding analysis indicates that, for a fixed system size \( L_z \) (including vacuum layers), the numerical error associated with neglecting the ELC term exhibits a non-trivial dependence on the number of image charge layers \( M \). 
Specifically, for \( |g_{\T u}g_{\T d}|>1 \), the error grows exponentially with \( M \);
for \( |g_{\T u}g_{\T d}| = 1 \), it grows linearly; and for \( |g_{\T u}g_{\T d}| < 1 \), no error escalation is observed as \( M \) increases.

\begin{figure}[htbp]
    \centering
    \includegraphics[width=0.55\linewidth]{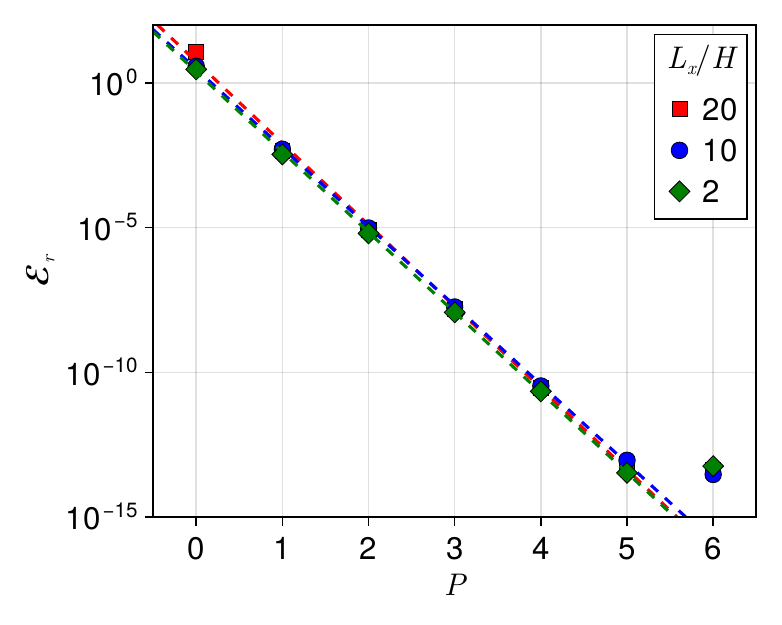}
    \caption{Relative errors in force ($\mathcal{E}_r$) as a function of the padding ratio $P$ for systems without dielectric interfaces. 
    We consider systems with heights $H = 0.5, 1, 5$ while fixing $L_x=L_y=10$. The padding ratio is defined as $P = (L_z - H) / L_x$. The dashed lines represent the fitted curves with decay rates using theoretical prediction Eq.~\eqref{eq::ForceError}. }
    \label{fig:elc_error_force}
\end{figure}

In what follows, we perform a series of numerical tests to validate our theoretical analysis. 
First, we examine the simplest case, i.e., systems \emph{without} dielectric interfaces by setting $\gamma_{\T u} = \gamma_{\T d} = 0$. 
In our numerical tests, we fix $L_x = L_y = 10$, and vary the system heights by setting $H = 0.5, 1, 5$. 
For clarity, we further introduce the dimensionless \emph{padding ratio} $P$, defined as $P = (L_z - H) / L_x$. 
The results, presented in \Cref{fig:elc_error_force}, clearly show that the relative errors in force $\mathcal{E}_r$ decay exponentially with $P$. 
Notably, the convergence with respect to $P$ is independent of the specific choice of $H$, aligning with our theoretical predictions in Eq.~\eqref{eq::ForceError}. 
Furthermore, our findings underscore the computational challenges of simulating \emph{strongly-confined} systems, characterized by a high aspect ratio, i.e., $L_x / H$. 
According to \Cref{fig:elc_error_force}, achieving single- or double-precision relative accuracy requires padding the system in the $z$ direction such that $(L_z - H)/L_x$ reaches around 3 and 5, respectively. For strongly-confined systems, this requires one to set $L_z\gg H$, which can significantly increases the computational cost if grid-based algorithms are used.

\begin{figure}[htbp]
    \centering
    \includegraphics[width=0.98\linewidth]{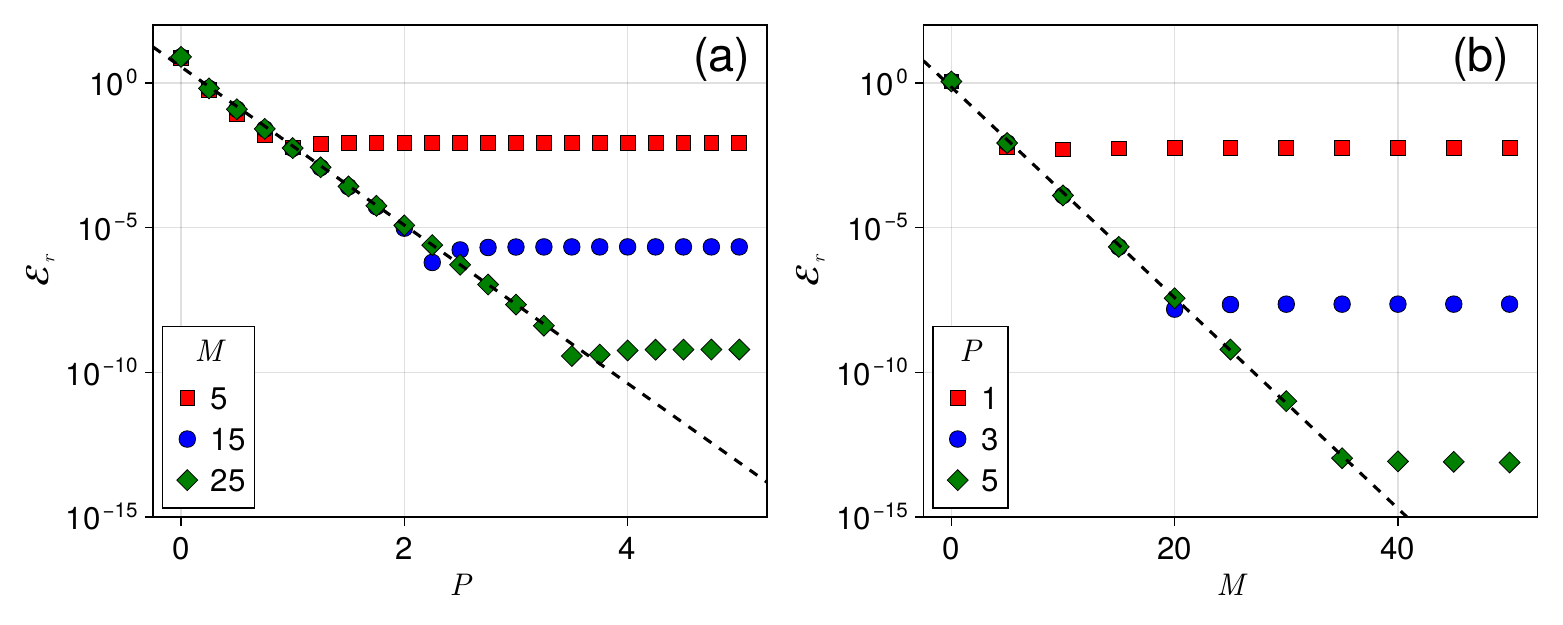}
    \caption{Relative errors in force ($\mathcal{E}_r$) for systems with dielectric interfaces. Here we fix $\gamma_{\T u}=\gamma_{\T d}=\gamma = 0.6$, $L_x=L_y=10$ and $H = 0.5$. 
    Panel (a) illustrates errors as a function of padding ratio $P$ with fixed image charge layers ($M=5,~15,~25$); panel (b) illustrates errors as a function of $M$ with fixed padding ratios ($P = 1,~3,~5$). The dashed lines in (a) and (b) represent the fitted curves with decay rates using theoretical predictions Eq.~\eqref{eq::lhs_less1} and Eq.~\eqref{eq::Ferr1}, respectively.}
    \label{fig:error_icm_pad_gamma_0.6_force}
\end{figure}

Next, we examine systems with dielectric interfaces.
We set $\gamma_{\T u}=\gamma_{\T d}=\gamma$ and explore two prototypical scenarios by choosing $\gamma = 0.6$ and $1$, respectively. 
In both scenarios, the size of simulation box is fixed as $L_x=L_y=10$ and $H = 0.5$, and we vary the padding ratio $P$ (by changing $L_z$) and the image series truncation parameter $M$ to validate our theoretical results.
We first address the scenario with $\gamma = 0.6$, where we have $|g_{\T u}g_{\T d}| < 1$, so that the numerical errors are mainly contributed by the image series truncation error Eq.~\eqref{eq::Ferr1} and the ELC term Eq.~\eqref{eq::lhs_less1}.
The numerical results, as shown in \Cref{fig:error_icm_pad_gamma_0.6_force}, demonstrate that the errors in force decay exponentially with both the padding ratio $P$ and the number of image charge layers $M$, consistent with our theoretical predictions. Notice that in \Cref{fig:error_icm_pad_gamma_0.6_force} (a), the errors saturate at certain accuracy levels, this is due to the fixed image series truncation error (as we fix $M$). Similarly, the errors saturate in \Cref{fig:error_icm_pad_gamma_0.6_force} (b) as they reach the fixed ELC errors for given padding ratios $P$.

\begin{figure}[htbp]
    \centering
    \includegraphics[width=0.98\linewidth]{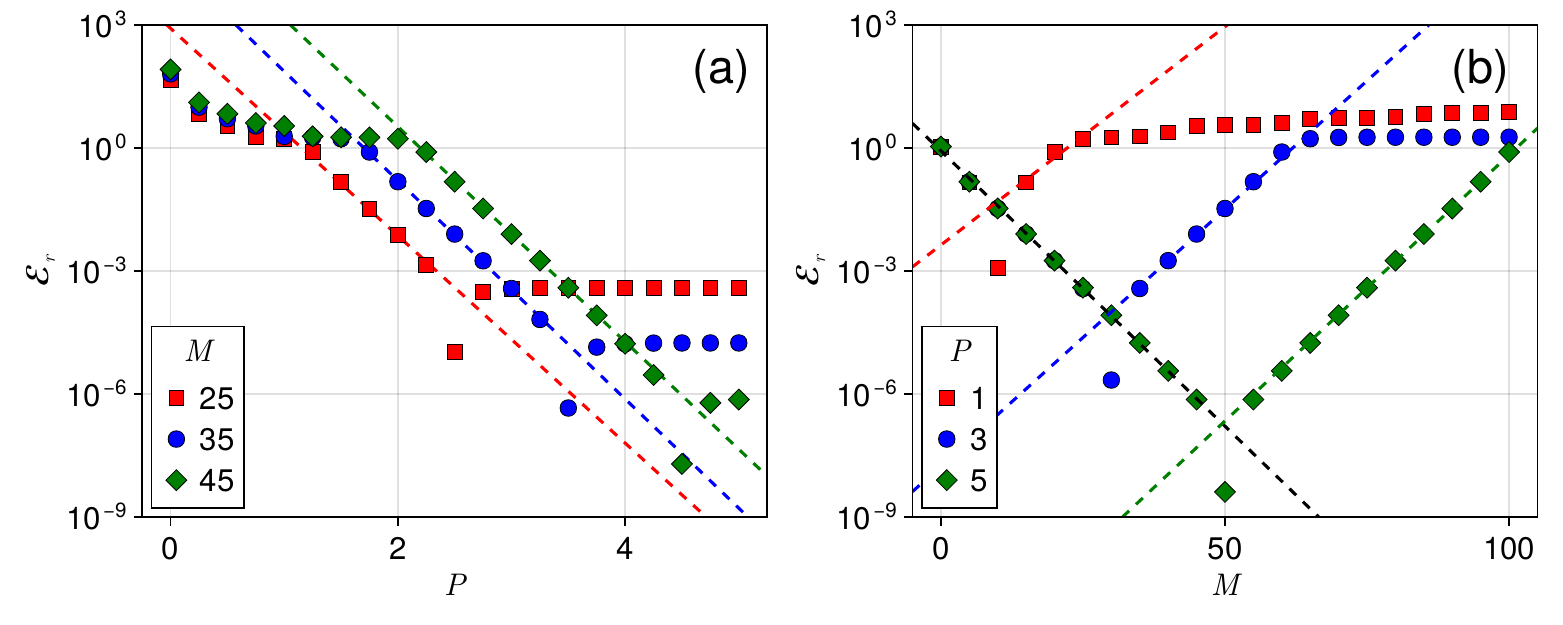}
    \caption{Relative errors in force ($\mathcal{E}_r$) for systems with dielectric interfaces. Here we fix $\gamma_{\T u}=\gamma_{\T d}=\gamma = 1$, $L_x=L_y=10$ and $H = 0.5$. 
    Panel (a) illustrates errors as a function of padding ratio $P$ with fixed image charge layers ($M=25,~35,~45$); panel (b) illustrates errors as a function of $M$ with fixed padding ratios ($P = 1,~3,~5$). The dashed lines in (a) and (b) represent the fitted curves with decay/growth rates using theoretical predictions Eq.~\eqref{eq::lhs_bigger1} and Eq.~\eqref{eq::Ferr1}.}
    \label{fig:error_icm_pad_gamma_1_force}
\end{figure}

Now consider the scenario with $\gamma = 1$, where we have $|g_{\T u}g_{\T d}| > 1$, so that the numerical errors associated with the ELC term are described by Eq.~\eqref{eq::lhs_bigger1}.
In this context, the error behavior becomes more complex: as $M$ increases, the ELC error exhibits exponential growth, whereas the image truncation error decreases exponentially.
As illustrated in \Cref{fig:error_icm_pad_gamma_1_force} (a), by fixing $M$, we still observe exponential convergence in forces as the padding ratio $P$ increases. 
However, for a fixed padding ratio $P$, the errors display non-monotonic behavior with increasing $M$, as depicted in \Cref{fig:error_icm_pad_gamma_1_force} (b).
This subtle phenomenon was not fully understood since its first observation in the work of Yuan \emph{et al.}~\cite{yuan2021particle}. 
Our analysis reveals it as the combined effect of image truncation and ELC errors: 
initially, errors decrease due to the decay of image truncation errors; and once $M$ surpasses a certain $P$-dependent threshold, the error amplification mechanism predicted by Eq.~\eqref{eq::lhs_bigger1} becomes dominant, causing errors to grow exponentially.
Finally, we note that in strongly-confined systems, the presence of dielectric interfaces would introduce additional computational challenges~\cite{dos2015electrolytes}. 
Specifically, for systems with higher aspect ratios, a larger $M$ is required to achieve the same level of accuracy. The relevant numerical results are summarized in Section 2 of the SI. 

\begin{figure}[htbp]
\centering
\includegraphics[width=0.98\linewidth]{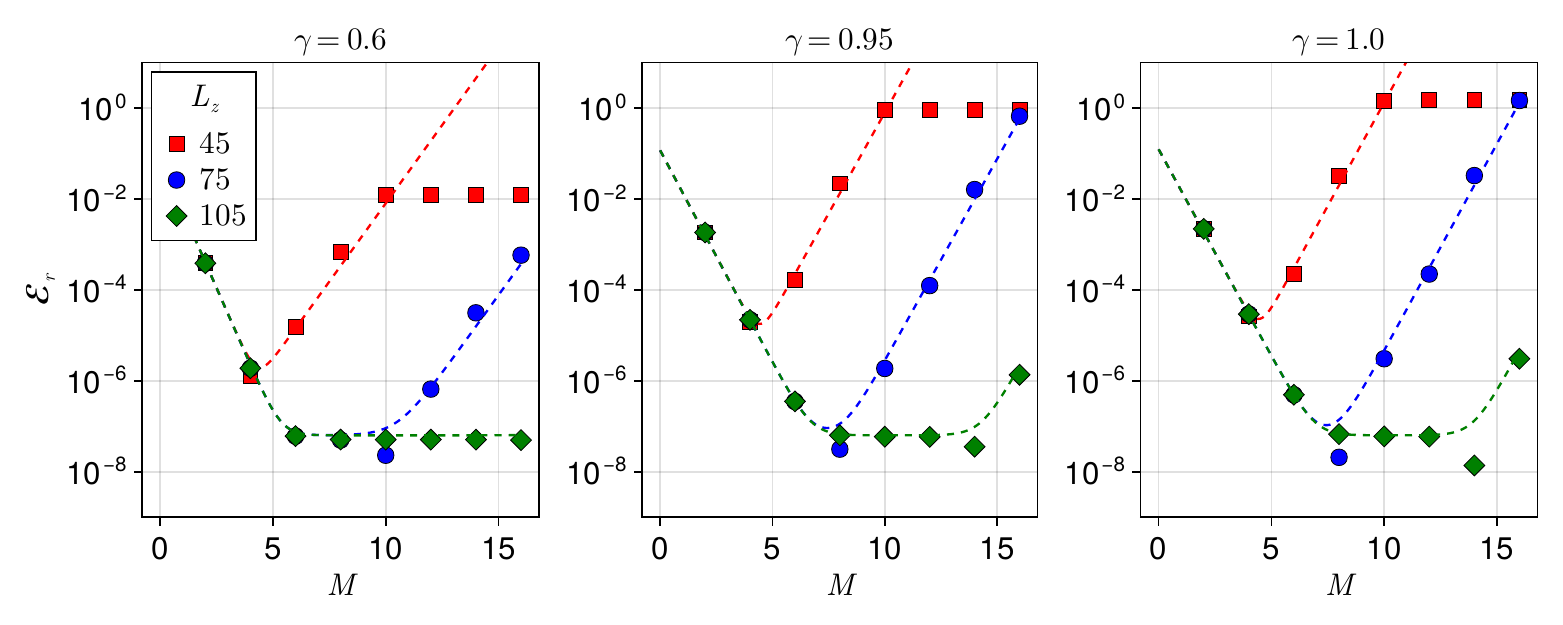}
\caption{
Relative errors ($\mathcal{E}_r$) in electrostatic energy for systems of 2:1 electrolytes with dieletric interfaces. 
Here we consider the same system setup studied by Yuan \emph{et al.}~\cite{yuan2021particle} using the ICM-PPPM method with $\gamma_{\T u}=\gamma_{\T d}=\gamma$ = 0.6, 0.95, and 1 (from left to right), respectively. In each panel, we fix $L_x = L_y = 15$, $H = 5$, and consider $L_z=$$45$, $75$, and $105$. 
Finally, the dashed lines represent the fitted curves using the sum of Eqs.~\eqref{eq:Uerr} and~\eqref{eq::U_ELC} (with coefficients determined by fitting).
}
\label{fig:error_yuan}
\end{figure}

As a final example, we validate our theoretical findings by replicating the non-monotonic error behavior reported by Yuan \emph{et al.} in figure 2 of Ref~\cite{yuan2021particle}. 
The reproduced numerical results are shown in \Cref{fig:error_yuan}, where we examine systems of dielectric-confined 2:1 electrolytes with $L_x = L_y = 15$ and $H = 5$. 
The three panels in \Cref{fig:error_yuan} correspond to systems with different reflection factors: $\gamma_{\T u} = \gamma_{\T d} = \gamma = 0.6$, $0.95$, and $1$, respectively. 
Within each panel, $L_z$ is varied as $45$, $75$, and $105$.
Note that for all cases considered here, the condition $|g_{\T u}g_{\T d}| > 1$ is satisfied, suggesting a similar phenomenon to that shown in \Cref{fig:error_icm_pad_gamma_1_force} (b), where errors diverge exponentially as $M$ surpasses a certain threshold.
Finally, we fit the numerical results using an analytical expression that combines the two primary error terms, Eq.~\eqref{eq:Uerr} and Eq.~\eqref{eq::U_ELC}, with coefficients determined through fitting.
The fitted curves, represented by dashed lines in \Cref{fig:error_yuan}, exhibit excellent agreement with the numerical data, thereby validating our theoretical predictions.

\section{Optimal parameter selection strategy}\label{sec:parameter}
In practical simulations of dielectric-confined systems, a systematic strategy for determining the optimal algorithm parameters is highly beneficial. Specifically, given a prescribed error tolerance $\varepsilon$, it is crucial to identify parameter choices that achieve this tolerance with minimum computational cost. 
For ICM-Ewald3D~\cite{dos2015electrolytes} and ICM-PPPM~\cite{yuan2021particle} algorithms, 
our error analysis indicates that four interrelated parameters must be determined to control accuracy: 1). the Ewald splitting parameter $\alpha$, 2). the real-space cutoff $r_c$, 3). the image charge series truncation parameter $M$, and 4). the padding length $L_z$. 
By combining the error estimates from this study with the established Ewald splitting error, the overall error estimates for the ICM-Ewald3D and ICM-PPPM methods can be expressed as follows:
\begin{equation}\label{eq:error_icmelc}
\resizebox{0.935\width}{!}{\(
\varepsilon \sim O\left(\frac{e^{-s^2}}{s^2}+\abs{\gamma_{\T u} \gamma_{\T d}}^{\lfloor\frac{M+1}{2}\rfloor} e^{-\frac{4\pi H \lfloor\frac{(M+1)}{2}\rfloor}{\max\{L_x,L_y\}}} + e^{-\frac{2\pi(L_z-H)}{\max\{L_x,L_y\}}} + \sum\limits_{l=1}^{M}C_{\gamma}^{(l)}e^{-\frac{2\pi(L_z-(l+1)H)}{\max\{L_x,L_y\}}}+e^{-\alpha^2 (L_z-H)^2}\right)\;,
\)}
\end{equation}
where the first term represents the Ewald decomposition error, the second term denotes the image charge series truncation error, the third and fourth terms correspond to the errors associated with the ELC term with image charges, and the last term corresponds to the trapezoidal discretization error.
Recall that, based on our analysis, the fourth term can grow exponentially with increasing $M$ if the condition $|g_{\T u}g_{\T d}| > 1$ is met. 
As a result, one should be careful in properly selecting the parameter $M$. 
Choosing a very large $M$ will not only increase the computational cost, but may also lead to incorrect results under certain circumstances.
In what follows, we propose an optimal parameter selection strategy based on the theoretical guidance of Eq.~\eqref{eq:error_icmelc}.

\emph{Step 1}: select $M$, so that the second term in Eq.~\eqref{eq:error_icmelc} is controlled by $\varepsilon$. 
There are three possible cases depending on the system setups:
\begin{enumerate}
    \item \textbf{Case 1}: If there is no dielectric interface, i.e., $\gamma_{\T u}=\gamma_{\T d}=0$, one simply set $M=0$.
    \item \textbf{Case 2}: If there is only one dielectric interface, i.e., either $\gamma_{\T u}=0$,$~\gamma_{\T d}\neq 0$ or $\gamma_{\T u}\neq 0$,$~\gamma_{\T d}=0$, one can set $M=1$ (since there is no image charge reflection).
    \item \textbf{Case 3}: If there are two polarizable dielectric interfaces, i.e., $\gamma_{\T u}\gamma_{\T d}\neq 0$, we select $M$ according to the following condition (obtained by algebraic manipulations of Eq.~\eqref{eq:Uerr}):
    \begin{equation}\label{eq::38}
    M\sim \frac{2\log \varepsilon - \frac{4\pi H}{\max\{L_x,L_y\}} - \log|\gamma_{\T u} \gamma_{\T d}|}{\log|\gamma_{\T u} \gamma_{\T d}| - \frac{4\pi H}{\max\{L_x,L_y\}}}\;.
\end{equation}
\end{enumerate}

\emph{Step 2}: select $L_z$, such that the sum of the third and fourth terms in Eq.~\eqref{eq:error_icmelc} is controlled by $\varepsilon$. 
Based on the leading-order error analysis presented in Section~\ref{sec::leadingerr}, two cases should be considered: 
\begin{enumerate}
    \item \textbf{Case 1}: If the condition $|g_{\T u}g_{\T d}|=\left|\gamma_{\T u}\gamma_{\T d}e^{\frac{4\pi H}{\max\{L_x,L_y\}}}\right|<1$
is satisfied, we have
\begin{equation}\label{eq::Lz_gamma_u_gamma_d_less1}
    L_z \geq H + \frac{\max\{L_x,L_y\}}{2\pi} \left( \log\frac{1}{\varepsilon} + \log \abs{\gamma_{\T u} + \gamma_{\T d} + e^{- \frac{2\pi H}{\max\{L_x,L_y\}}}} \right)\;.
\end{equation}
Note that for this case, $L_z$ can be chosen independent of $M$.
\item \textbf{Case 2}: If $\left|\gamma_{\T u}\gamma_{\T d}e^{\frac{4\pi H}{\max\{L_x,L_y\}}}\right|\geq 1$, we obtain
\begin{equation}\label{eq::Lz_gamma_u_gamma_d_greater1}
    L_z \geq (M + 1) H + \frac{\max\{L_x,L_y\}}{2\pi} \left( \log\frac{1}{\varepsilon} + \log \abs{\gamma_{\T u} \gamma_{\T d}} \right)\;.
\end{equation}
\end{enumerate}
It is interesting to note that, the selection of $L_z$ as derived in Eq.~\eqref{eq::Lz_gamma_u_gamma_d_greater1} can be interpreted physically as ensuring a sufficiently large vacuum layer in $z$ such that all the image charges can not overlap each other due to the periodic boundary conditions (necessitating $L_z\geq (M+1)H$). Additionally, an extra buffer zone is required, the length of which is determined by the specific tolerance $\varepsilon$.

\emph{Step 3}. After $M$ and $L_z$ are determined, the Ewald splitting parameter $\alpha$ and real-space cutoff $r_c$ can be selected. 
As indicated by Eq.~\eqref{eq:error_icmelc}, the Ewald decomposition error is independent of both the image truncation and the ELC error terms. 
The only extra constraint comes from the trapezoidal discretization error, which is typically minor due to its rapid decay with increasing $L_z$.
Consequently, the standard strategy can be employed to choose these parameters by setting $\varepsilon = e^{-s^2}/s^2$ and solving for $s$, where $s=r_c/\alpha$. 
The specific values of $r_c$ and $\alpha$ can then be adjusted to balance the computational costs between real-space and reciprocal-space calculations~\cite{frenkel2023understanding}.
Finally, to guarantee that the trapezoidal discretization error has been controlled, it is necessary to verify that the following condition is satisfied: 
\begin{equation}
    \alpha \geq (L_z-H)^{-1}\sqrt{\log\varepsilon^{-1}}\;.
\end{equation}
If this condition is not met, then $\alpha$ must be relaxed to fulfill this constraint.

\begin{figure}[!htbp]
    \centering
    \includegraphics[width=0.92\linewidth]{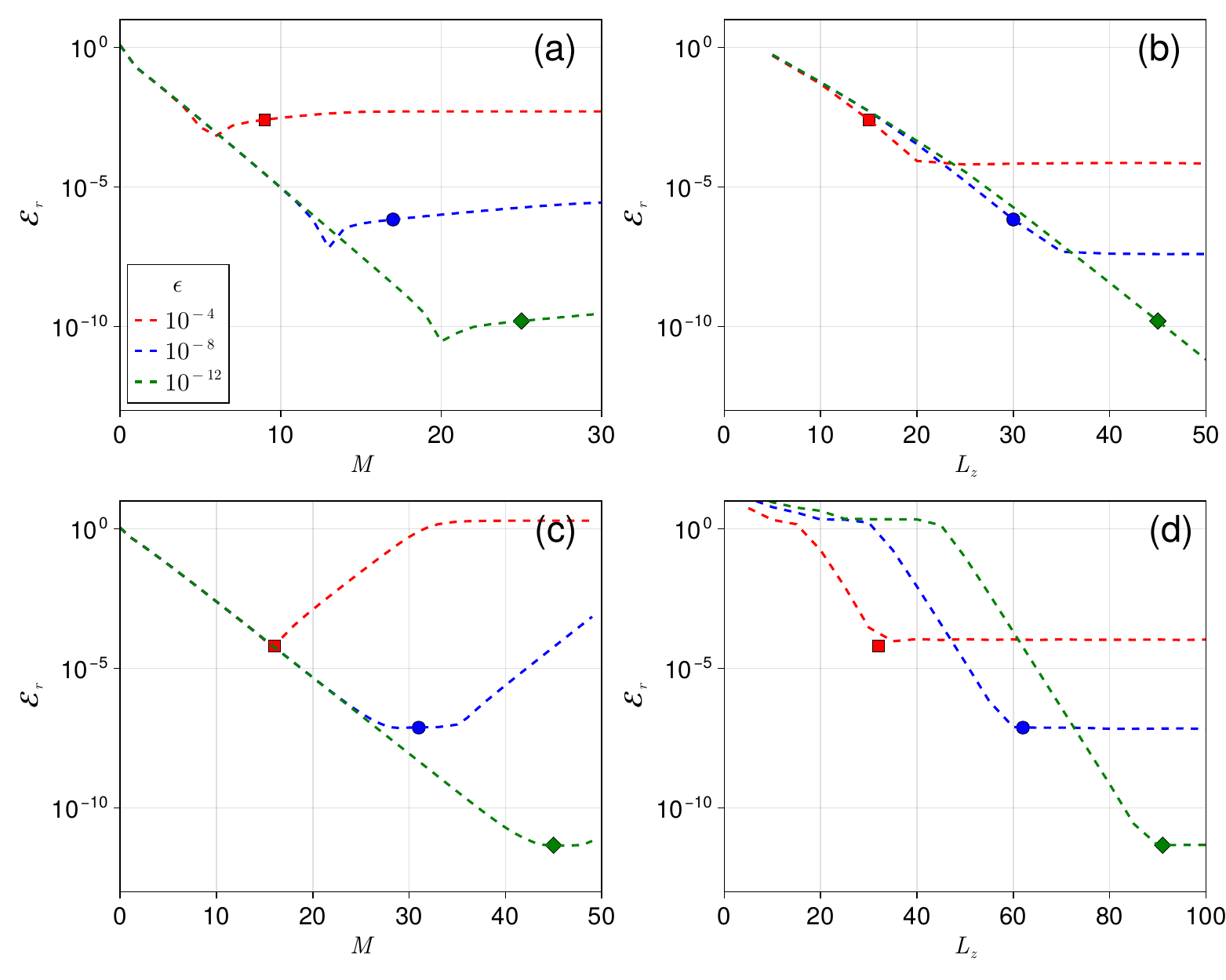}
    \caption{
        Relative errors in force ($\mathcal{E}_r$) for two prototypical dielectric-confined systems with $\gamma=0.6$ (panels (a-b)) and $\gamma = 1$ (panels (c-d)), respectively.
        For both systems, we fix $L_x=L_y=10$ and $H = 1$. 
        Within each panel, the dashed lines correspond to the numerical errors obtained according to tolerance values $\epsilon = 10^{-4}$, $10^{-8}$, and $10^{-12}$, and the solid markers indicate the specific choice of $M$ or $L_z$ selected via the proposed strategy.
        Note that in panels (a) and (c), $s$ and $L_z$ are fixed according to \Cref{tab:parameter_selection_results}, whereas in panels (b) and (d), we fix $s$ and $M$.
    }
    \label{fig:error_parameter_selection_force}
\end{figure}

We finally validate the proposed parameter selection strategy through numerical tests on two prototypical dielectric-confined systems, characterized by $\gamma_{\T u} = \gamma_{\T d} = 0.6$ and 1, respectively. 
The system dimensions are fixed at $L_x=L_y=10$ and $H = 1$. 
Note that for the system with $\gamma=0.6$, we have $|g_{\T u}g_{\T d}|<1$, while for the system with $\gamma=1$, $|g_{\T u}g_{\T d}|>1$.
For each system, the error tolerance $\varepsilon$ is set to $\epsilon = 10^{-4}$, $10^{-8}$, and $10^{-12}$. 
By applying the parameter selection strategy discussed earlier, we are able to select the optimized parameters and the results are summarized in \Cref{tab:parameter_selection_results}. 
Detailed error curves as functions of $M$ and $L_z$ are plotted in \Cref{fig:error_parameter_selection_force}, where the solid markers denote the specific parameter chosen via the proposed strategy.
\begin{table}[htbp]
    \centering
    \begin{tabular}{|c|c|c|c|c|}
        \hline
        $\gamma$ & $\epsilon$ & $s$ & $M$ & $L_z$ \\
        \hline
        $0.6$ & $10^{-4}$ & 3 & 9 & 15 \\
        $0.6$ & $10^{-8}$ & 4 & 17 & 30 \\
        $0.6$ & $10^{-12}$ & 5 & 25 & 45 \\
        $1$ & $10^{-4}$ & 3 & 16 & 32 \\
        $1$ & $10^{-8}$ & 4 & 31 & 62 \\
        $1$ & $10^{-12}$ & 5 & 45 & 91 \\
        \hline
    \end{tabular}
    \caption{Algorithm parameters $s$,~$M$ and $L_z$ for varied tolerance $\epsilon = 10^{-4}$, $10^{-8}$, and $10^{-12}$. 
    The parameters are selected according to the strategy proposed in this work. Two prototypical dielectric-confined systems are considered, with $\gamma=0.6$ and 1, respectively. Both systems have dimensions $L_x=L_y=10$ and $H = 1$.} \label{tab:parameter_selection_results}
\end{table}
It is evident from \Cref{fig:error_parameter_selection_force} that, across all test cases with varying $\gamma$ and $\varepsilon$, the selected parameters consistently achieve optimal or near-optimal performance, especially for the case $\gamma = 1$. 
Consequently, we conclude that the parameter selection strategy proposed herein offers practical guidance for optimizing the performance of MD simulations of dielectric-confined systems.

\section{Conclusion}\label{sec:conclusion}
In this work, we presented a rigorous error analysis of Ewald summation for dielectric-confined planar systems, where the polarization potential and force field are modeled using an infinitely reflected image charge series. 
In particular, we address the truncation error of the image charge series and the error estimations associated with the ELC term involving image charges, which may introduce significant errors but are often over-looked.
Our error estimations are validated numerically across several prototypical systems. Moreover, through analysis, we are able to elucidate the counterintuitive non-monotonic error behavior observed in previous simulation studies.
Finally, based on the theoretical insights, we propose an optimal parameter selection strategy, offering practical guidance for achieving efficient and accurate MD simulations of dielectric-confined systems.





\section*{Acknowledgement}\label{sec:acknowledgement}
The work of X. G. and Z. G. is supported by the Natural Science Foundation of China (Grant No. 12201146), Natural Science Foundation of Guangdong (Grant No. 2023A1515012197), and the Guangzhou-HKUST(GZ) Joint Research Project (Grant Nos. 2023A03J0003 and 2024A03J0606). 
The work of Q. Z. and J. L. are supported by the Natural Science Foundation of China (grant Nos. 12325113, 12426304 and 12401570) and the Science and Technology Commission of Shanghai Municipality (grant No. 23JC1402300). The work of J. L. is partially supported by the China Postdoctoral Science Foundation (grant No. 2024M751948). The authors also acknowledge the support from the HPC center of Shanghai Jiao Tong University.



\section*{Data Availability}\label{sec:data}

Our software package \texttt{EwaldSummations.jl}~\cite{EwaldSummations} is developed based on the Julia Programming Language~\cite{Julia-2017} and the open-source package \texttt{CellListMap.jl}~\cite{celllistmap}, and the figures are generated by the open-source package \texttt{Makie.jl}~\cite{DanischKrumbiegel2021}.
The code and data used in this work are available at \url{https://github.com/ArrogantGao/ICMErrorEstimation}.

\appendix

\section{Numerical validations for error estimations in energy}\label{sec:numeric_energy}

\begin{figure}[htbp]
    \centering
    \includegraphics[width=0.98\linewidth]{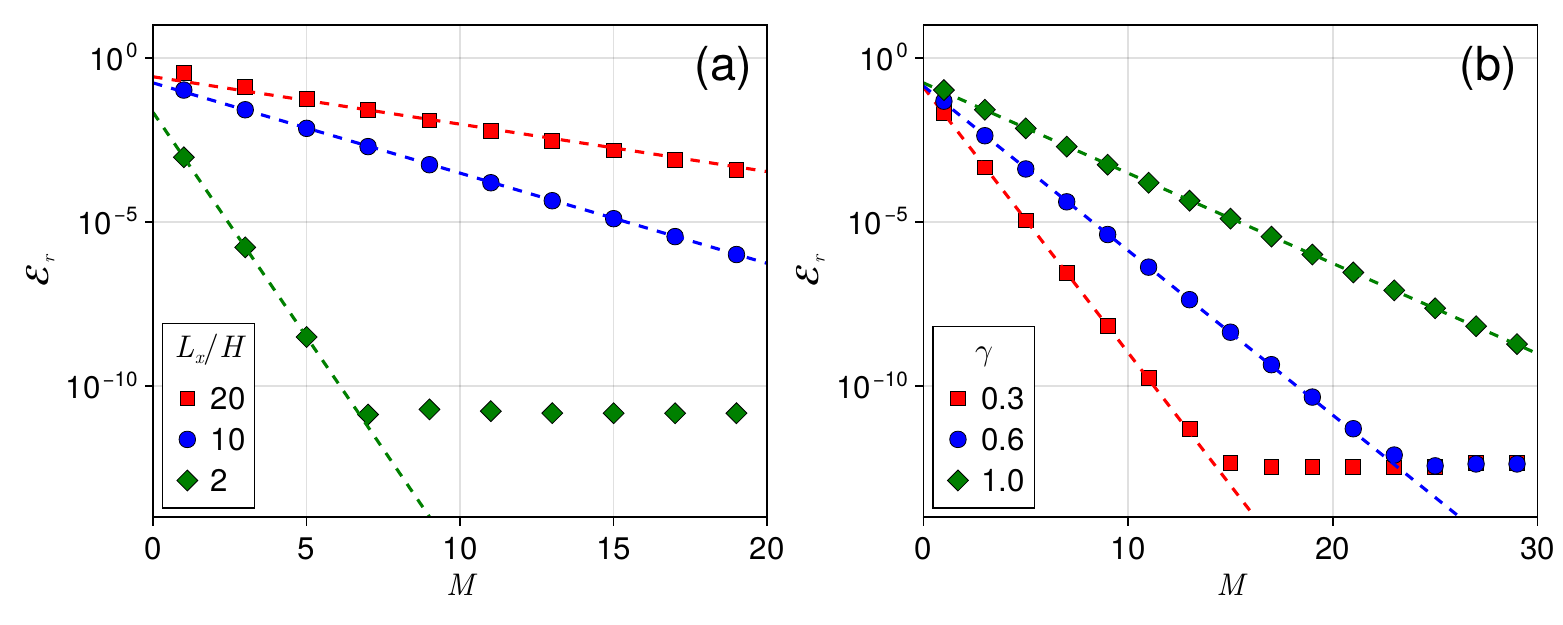}
    \caption{
        Relative errors of energy due to truncation of image charges are presented. The dashed lines represent the fitted decay rate, as described in Eq. (54) in the main text. In panel (a) we set $\gamma_{\T u} = \gamma_{\T d} = 1$ and consider system heights of $H = 0.5, 1$ and $5$. In panel (b), we set $H = 1$ while varying $\gamma_{\T u} = \gamma_{\T d} = \gamma = 0.3, 0.6$ and $1$. 
    }
    \label{fig:icm_error}
\end{figure}

\begin{figure}[htbp]
    \centering
    \includegraphics[width=0.49\linewidth]{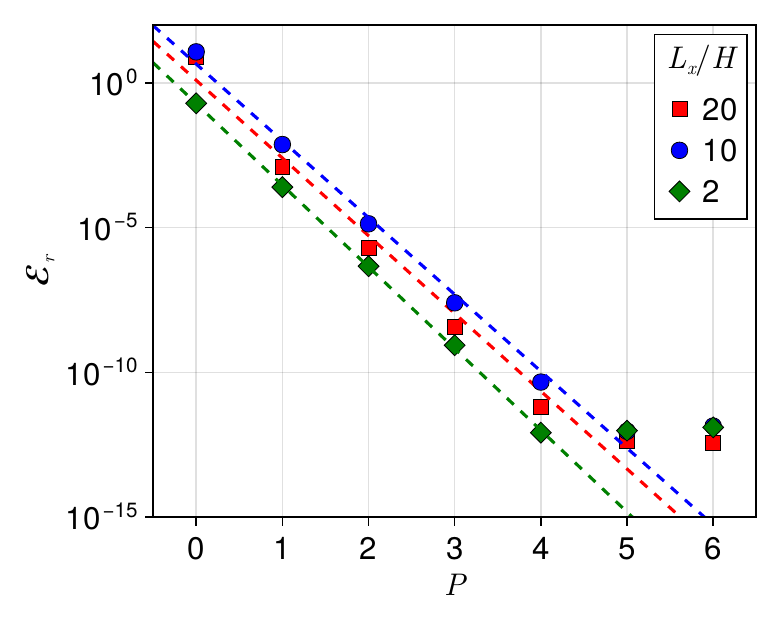}
    \caption{Relative error of energy of reformulating the Ewald2D summation as 3D ones are presented. We set $\gamma_{\T u} = \gamma_{\T d} = 0$ and consider system heights of $H = 0.5, 1, 5$. Here $P = (L_z - H) / L_x$ denotes the padding ratio.}
    \label{fig:elc_error}
\end{figure}

\begin{figure}[htbp]
    \centering
    \includegraphics[width=0.98\linewidth]{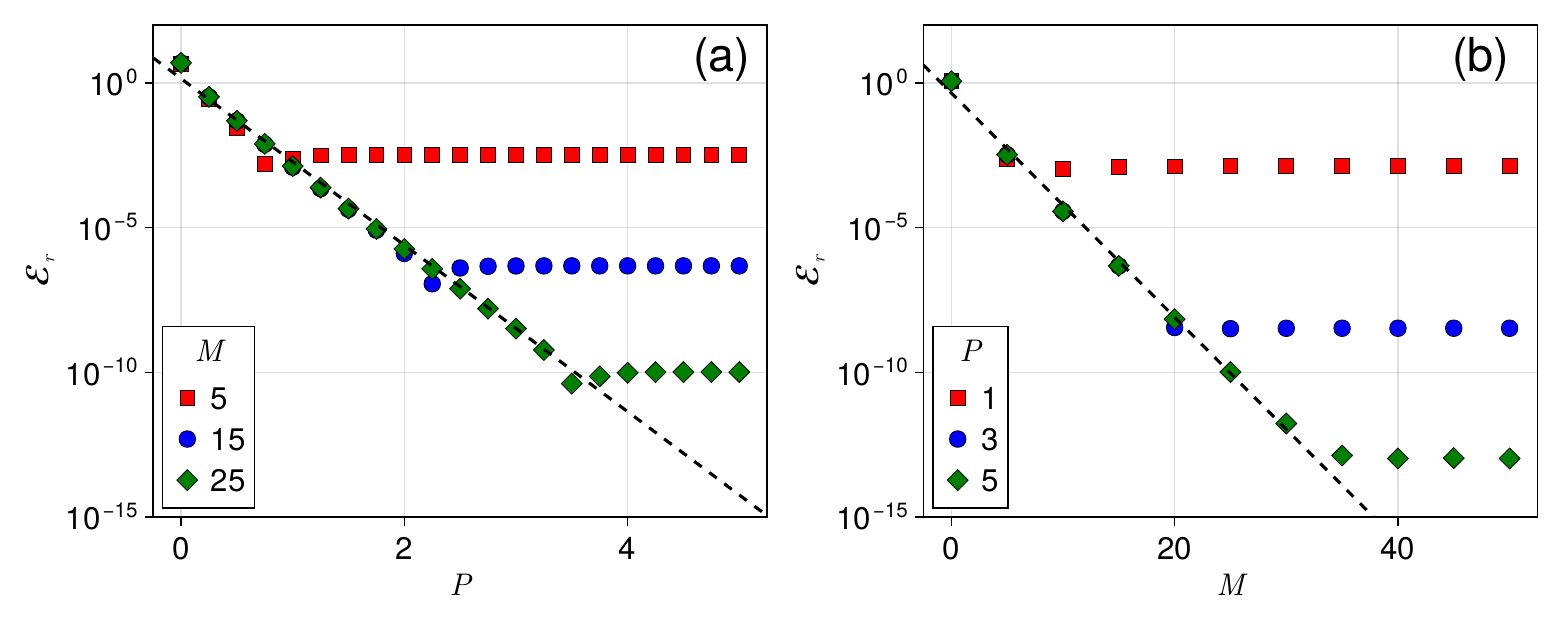}
    \caption{Relative error of energy in dielectric-confined Coulomb system with parameters $\gamma_{\T u} = \gamma_{\T d} = \gamma = 0.6$ and $H = 0.5$. Panel (a) illustrates error evolution with fixed image charge layers ($M = 5,~15,~25$) under varying padding ratios ($P$), whereas panel (b) demonstrates error progression with fixed $P = 1,~3,~5$ across increasing $M$.}
    \label{fig:error_icm_pad_gamma_0.6}
\end{figure}

\begin{figure}[htbp]
    \centering
    \includegraphics[width=0.98\linewidth]{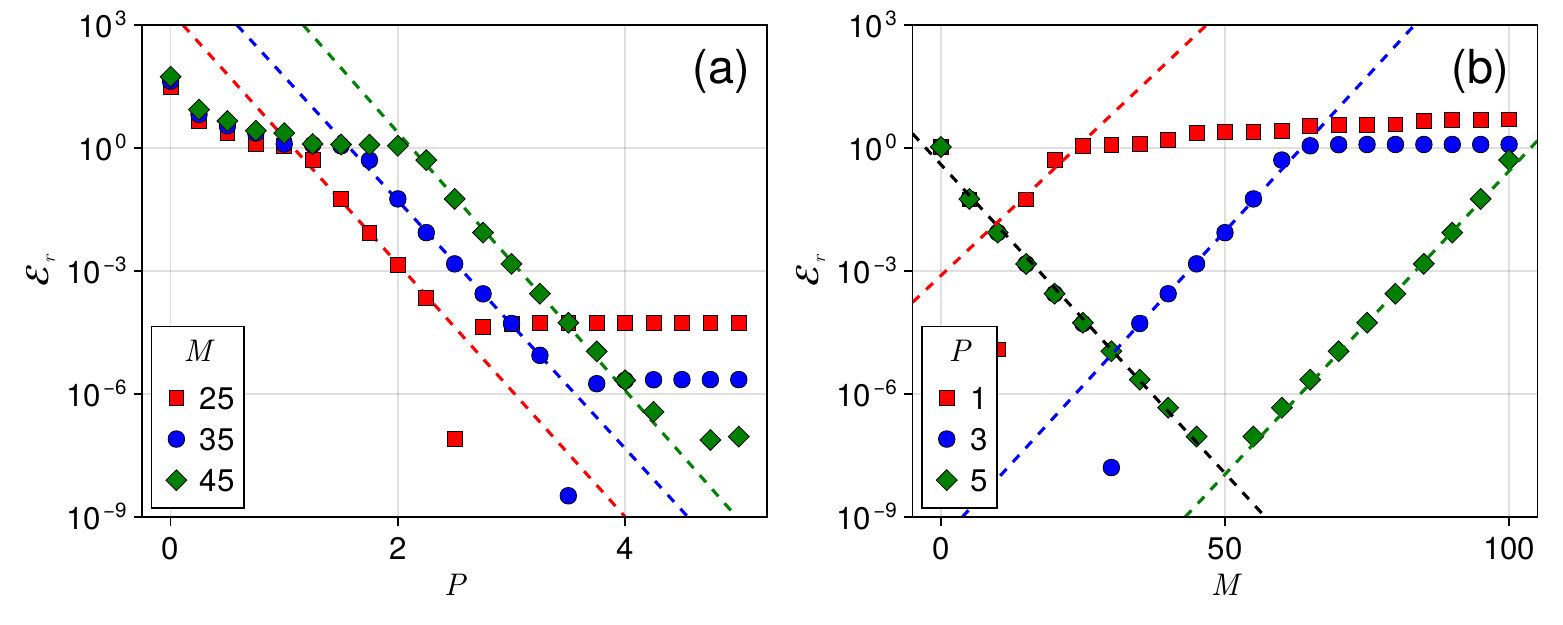}
    \caption{Relative error of energy in dielectric-confined Coulomb system with parameters $\gamma_{\T u} = \gamma_{\T d} = \gamma = 1$ and $H = 0.5$. Panel (a) illustrates error evolution with fixed image charge layers ($M = 5,~15,~25$) under varying padding ratios ($P$), whereas panel (b) demonstrates error progression with fixed $P = 1,~3,~5$ across increasing $M$.}
    \label{fig:error_icm_pad_gamma_1}
\end{figure}

In this section, we present supplementary results for the numerical validations for error estimations in electrostatic energy, complementing the main text. Figures~\ref{fig:icm_error} to~\ref{fig:error_icm_pad_gamma_1} correspond one-to-one with Figures 1 to 4 in the main text, maintaining the same systems and parameter settings as those used in the force calculation results. The dashed lines in these figures represent the fitted curves based on our theoretical estimation given in Eq.~(52) of the main text. The energy error results exhibit similar trends to the force errors and align well with our theoretical predictions.

\section{Challenge associated with strongly-confined systems}
\begin{figure}[htbp]
    \centering
    \includegraphics[width=0.98\linewidth]{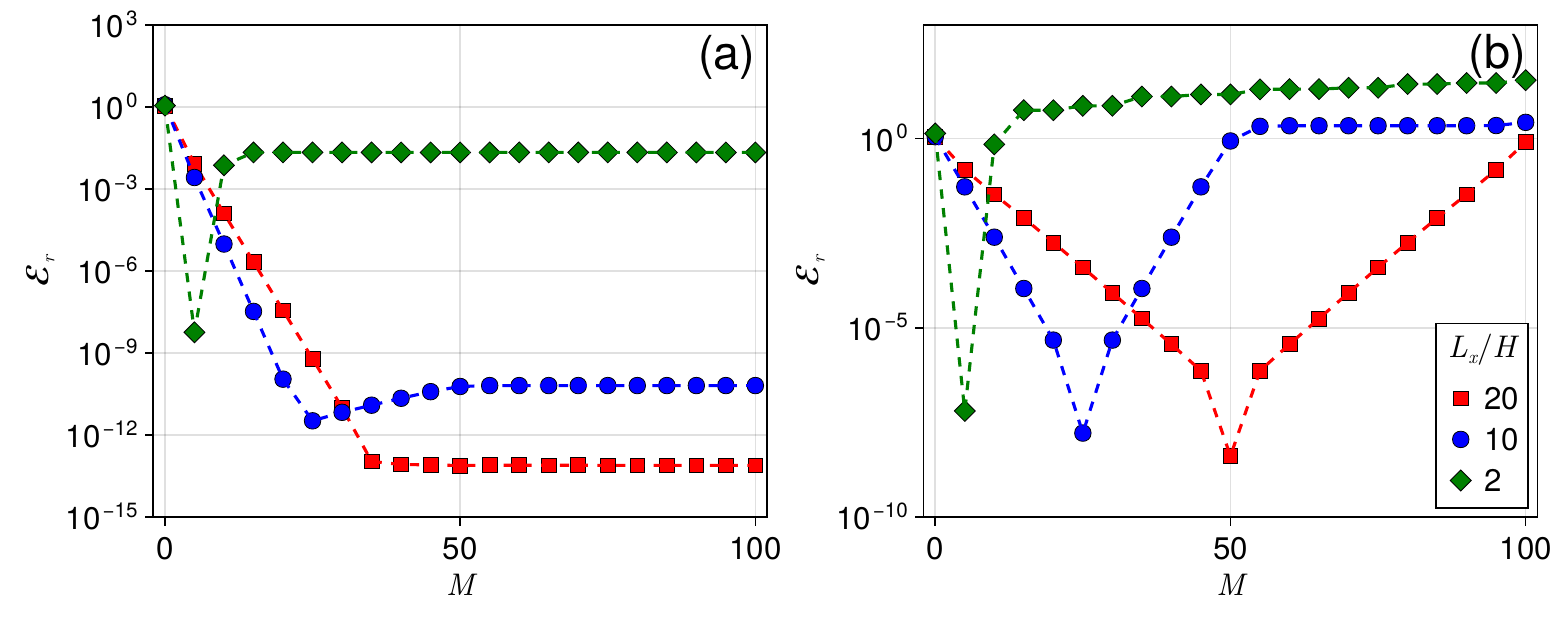}
    \caption{Relative error of force ($\mathcal{E}_r$) in a dielectric-confined Coulomb system with parameters $P = 5$ and $H = 0.5, 1, 5$. In panels (a) and (b), the dielectric contrasts are set to $\gamma_{\T u} = \gamma_{\T d} = \gamma = 0.6$ and $\gamma_{\T u} = \gamma_{\T d} = \gamma = 1$, respectively.}
    \label{fig:icm_elc_error_force}
\end{figure}

\begin{figure}[htbp]
    \centering
    \includegraphics[width=0.98\linewidth]{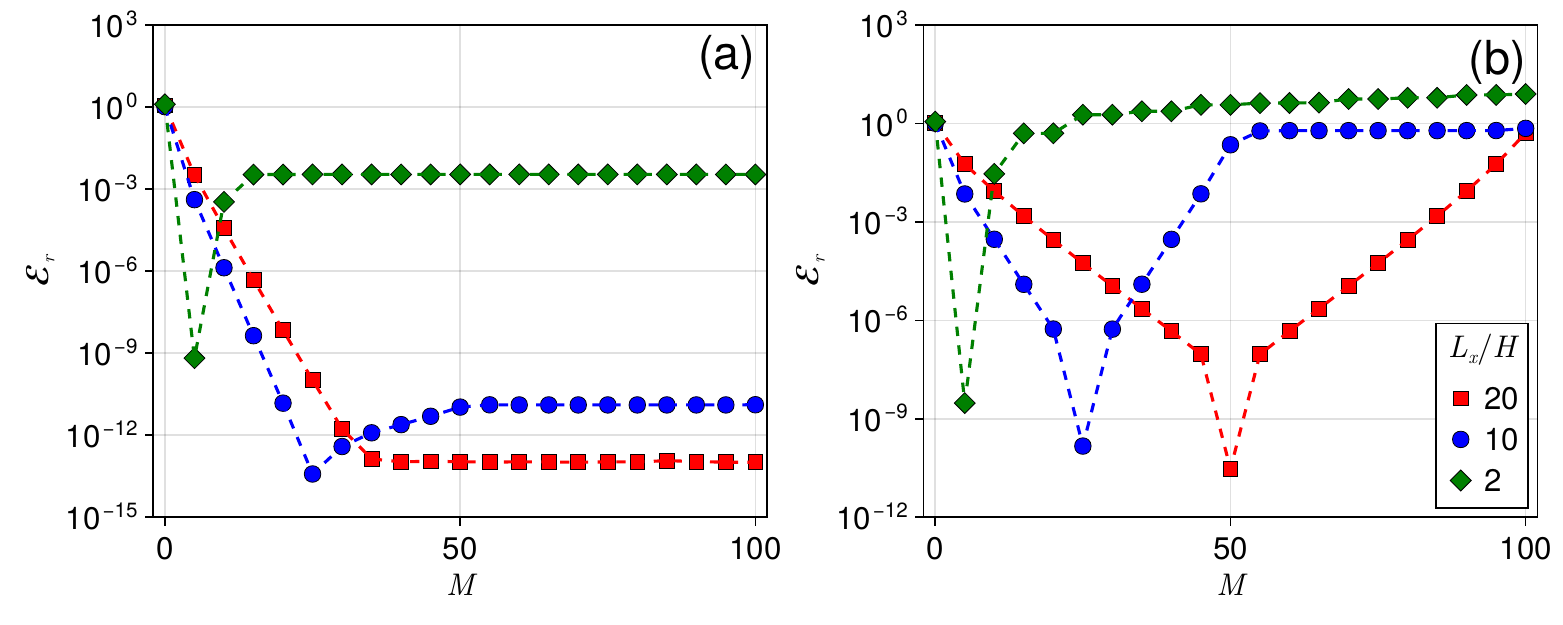}
    \caption{Relative error of energy in a dielectric-confined Coulomb system with parameters $P = 5$ and $H = 0.5, 1, 5$. In panels (a) and (b), the dielectric contrasts are set to $\gamma_{\T u} = \gamma_{\T d} = \gamma = 0.6$ and $\gamma_{\T u} = \gamma_{\T d} = \gamma = 1$, respectively.}
    \label{fig:icm_elc_error}
\end{figure}

In practical simulations, such as when studying thin membranes, ion transport in slit channels and supercapacitors, accurately capturing the effects of nanoconfinement, i.e., when $L_{x,y} / H\gg 1$, is crucial.
Previous numerical studies have found that more image layers are needed to achieve satisfactory accuracy for confined systems~\cite{dos2015electrolytes}. 
To further investigate the numerical properties of strongly-confined systems, we present the errors in force in \Cref{fig:icm_elc_error_force}, where we fix $L_x=L_y=10$, $P=(L_z-H)/L_x = 5$ and consider system heights $H = 0.5, 1, 5$ while varying the number of image charge layers $M$. 
In \Cref{fig:icm_elc_error_force} (a), for $\gamma = 0.6$, we observe that the error decays exponentially as $M$ increases for $H = 0.5$ and $H = 1$. 
However, for $H = 5$, where $|g_{\T u}g_{\T d}|>1$, the error becomes non-monotonic as $M$ increases, 
which is consistent with our theoretical predictions as discussed in the main text.
In \Cref{fig:icm_elc_error_force} (b), with $\gamma = 1$, we observe a similar non-monotonic pattern for all aspect ratios. 
It is important to note that, the rate of error decrease/increase depends on the aspect ratio $L_x/H$. 
A higher aspect ratio leads to a slower increase or decrease in errors, as $M$ is less or greater than the critical value that minimizes the error, respectively. 
This also highlights the computational challenge associated with strongly-confined systems -- to achieve the same accuracy, a much larger $M$ is required compared to non-confined systems.
Same conclusions hold for the errors in energy, as shown in \Cref{fig:icm_elc_error}.

\bibliography{reference}

\end{document}